\documentclass[11pt,reqno]{amsart}
\usepackage[T1]{fontenc}
\usepackage{lmodern, comment}
\usepackage{amsmath,amssymb,amsthm,mathtools}
\usepackage[expansion=false]{microtype}
\usepackage{enumitem}
\usepackage[sort,nocompress,noadjust]{cite}
\usepackage[hidelinks]{hyperref}
\usepackage{aliascnt}
\usepackage[nameinlink,noabbrev]{cleveref}
\usepackage[letterpaper,margin=1in]{geometry}
\numberwithin{equation}{section}
\newtheorem{theorem}{Theorem}[section]
\crefname{theorem}{Theorem}{Theorems}
\newaliascnt{lemma}{theorem}
\newtheorem{lemma}[lemma]{Lemma}
\aliascntresetthe{lemma}
\crefname{lemma}{Lemma}{Lemmas}
\newaliascnt{proposition}{theorem}
\newtheorem{proposition}[proposition]{Proposition}
\aliascntresetthe{proposition}
\crefname{proposition}{Proposition}{Propositions}
\newaliascnt{corollary}{theorem}
\newtheorem{corollary}[corollary]{Corollary}
\aliascntresetthe{corollary}
\crefname{corollary}{Corollary}{Corollaries}
\newaliascnt{conjecture}{theorem}

\aliascntresetthe{conjecture}
\crefname{conjecture}{Conjecture}{Conjectures}
\theoremstyle{definition}
\newaliascnt{example}{theorem}
\newtheorem{example}[example]{Example}
\aliascntresetthe{example}
\crefname{example}{Example}{Examples}
\newaliascnt{definition}{theorem}

\aliascntresetthe{definition}
\crefname{definition}{Definition}{Definitions}
\theoremstyle{remark}
\newtheorem*{remark*}{Remark}
\newaliascnt{remark}{theorem}
\newtheorem{remark}[remark]{Remark}
\aliascntresetthe{remark}
\crefname{remark}{Remark}{Remarks}
\newcommand{\N}{\mathbb N}
\newcommand{\Nzero}{\mathbb N_0}
\newcommand{\R}{\mathbb R}
\newcommand{\C}{\mathbb C}
\newcommand{\T}{\mathbb T}
\newcommand{\Z}{\mathbb Z}
\DeclareMathOperator{\HS}{HS}
\newcommand{\HH}{\mathcal H}
\newcommand{\YY}{\mathcal Y}
\newcommand{\BB}{\mathcal B}
\newcommand{\cA}{\mathcal A}
\newcommand{\ran}{\operatorname{ran}}
\newcommand{\rank}{\operatorname{rank}}

\newcommand{\norm}[1]{\left\lVert#1\right\rVert}
\newcommand{\abs}[1]{\left\lvert#1\right\rvert}
\newcommand{\ip}[2]{\left\langle#1,#2\right\rangle}
\newcommand{\Dens}{\mathfrak d}
\DeclareMathOperator{\diag}{diag}
\title[Totally positive interpolation and critical sampling]{Totally Positive Interpolation and Critical Sampling\\of Positive Operator Orbits}
\author[I. A. Krishtal]{Ilya A. Krishtal}
\address{School of Mathematical and Statistical Sciences, Northern Illinois University, DeKalb, Illinois 60115, USA}
\email{ikrishtal@niu.edu}
\author[J. Wu]{Jian Wu}
\thanks{Corresponding author: Jian Wu (\texttt{xingxingwu2022@163.com}).}
\address{School of Mathematics, Southwestern University of Finance and Economics, Chengdu 611130, People's Republic of China}
\email{xingxingwu2022@163.com}
\date{}
\hypersetup{pdftitle={Totally Positive Interpolation and Critical Sampling of Positive Operator Orbits},pdfauthor={Ilya A. Krishtal and Jian Wu}}

\begin{document}
\begin{abstract}
How few observations are needed for the stable recovery of the initial state of a linear dynamical system? We study this question for systems governed by a positive operator, with fixed sensors and freely chosen sampling times and positive weights.
Our main tool is an interpolation theorem for translates of a nonzero continuous integrable totally positive function. For every uniformly separated sequence of centers, we choose one node near each center so that the corresponding collocation operator is boundedly invertible. The displacement, node separation, and operator bounds depend only on the kernel and the lower separation of the centers; no upper bound is needed on their gaps. If a single integer orbit is a frame, we construct weighted Riesz bases whose sampling counts attain
the spectral lower bound up to a bounded additive error.
The same method yields bases from continuous functional calculus, including explicit resolvent bases. In the case of normal operators, we give a single-orbit frame for which every admissible weighted design has sampling count growing faster than the radial spectral count.
\end{abstract}
\maketitle
\raggedbottom

\section{Introduction}\label{sec:intro}
We consider a bounded positive  operator $A$ on a complex separable Hilbert space $\HH$ and a bounded linear map $C:\HH\to\YY$ into another complex Hilbert space $\YY$. We regard $A$ as the evolution operator and $C$ as the observation map. The pair $(A,C)$ is \emph{observable} if
\begin{equation}\label{eq:original-frame}
 a\norm{x}^{2}\le\sum_{n=0}^{\infty}\norm{CA^n x}^{2}\le b\norm{x}^{2},
 \qquad x\in\HH,
\end{equation}
for some $0<a\le b<\infty$. The lower bound guarantees stable recovery of the initial state, while the upper bound controls the total energy of the observations.

For finite-dimensional output, observability is equivalent to a frame condition. Indeed, when $\YY=\C^r$, we write $Cx=(\ip{x}{f_\ell})_{\ell=1}^r$. Since $A$ is selfadjoint, the identity
\[
 \ip{A^n x}{f_\ell}=\ip{x}{A^n f_\ell}
\]
moves the evolution from the unknown state to the observation vectors. Thus \eqref{eq:original-frame} says precisely that the multiorbit $(A^n f_\ell)_{n\ge0,\,1\le\ell\le r}$ is a frame, as in \cite{ACCMP17,ACMT17,CMPP20}. This recovery problem also appears in control theory as exact observability of diagonal systems \cite{JZ01}. Continuous powers and the relation between continuous and discrete dynamical sampling are studied in \cite{AHP19,DMM21} and reviewed in the survey \cite{ACKM26}.

Here we ask how many observation times are needed when $A$ and $C$ remain fixed. We choose the times and their positive weights, recording both in a sampling design
\[
 \mu=\sum_{t\in T}w_t\delta_t,
 \qquad T\subset[0,\infty)\text{ locally finite},\quad w_t>0.
\]
Here $\delta_t$ is the unit point mass at $t$, and local finiteness means that $T\cap[0,R]$ is finite for every $R\ge0$. The powers $A^t$, $t>0$, are defined by the spectral calculus, with $A^0=I$ on all of $\HH$. At each time, one weight is used for all output coordinates. The counting function and the formal sampled Gramian are
\[
 N_\mu(R)=\#(T\cap[0,R])\quad\mbox{and}\quad
 G_\mu=\sum_{t\in T}w_t A^tC^*CA^t.
\]
The sum is understood through the nonnegative quadratic form $q_\mu(x)=\sum_{t\in T}w_t\norm{CA^tx}^2$, which may be infinite. If $q_\mu(x)\le\beta\norm x^2$, the analysis map $U_\mu x=(\sqrt{w_t}\,CA^tx)_{t\in T}$ is bounded into $\ell^2(T;\YY)$ and $G_\mu=U_\mu^*U_\mu$, with strong convergence of the partial sums. We use this convention throughout.

A design is \emph{admissible} if $\alpha I\le G_\mu\le\beta I$ for some $0<\alpha\le\beta<\infty$. We denote the admissible real and nonnegative-integer designs by $\cA_\R(A,C)$ and $\cA_\Z(A,C)$, respectively. Each time contributes one to $N_\mu$, regardless of the number of scalar observations made at that time. When only the support matters, we also write $N_T(R)=\#(T\cap[0,R])$.

To compare designs at the time scale $R$, we count the spectral modes with decay rates at least $1/R$. If $E_A$ is the projection-valued spectral measure of $A$, this count is
\begin{equation}\label{eq:spectral-count}
 n_A(R)=\rank E_A((0,e^{-1/R}]),\qquad R>0.
\end{equation}
Whenever $n_A(R)<\infty$ for every $R>0$ and $n_A(R)\to\infty$ as $R\to\infty$, we define
\[
 c(\mu)=\limsup_{R\to\infty}\frac{N_\mu(R)}{n_A(R)}\quad \mbox{and}
\quad\Dens_\ast(A,C)=\inf_{\mu\in\cA_\ast(A,C)}c(\mu),\quad \ast\in\{\R,\Z\}.
\]
The zero eigenspace is not included in $n_A$. The infimum $\Dens_\ast(A,C)$ is \emph{attained} if $c(\mu)=\Dens_\ast(A,C)$ for some admissible design. The stronger conclusion
\[
 N_\mu(R)=\Dens_\ast(A,C)n_A(R)+O(1)
\]
is called \emph{critical-count attainment}. Throughout this paper, $O(1)$ denotes an error bounded as $R\to\infty$; subscripts indicate parameters on which the bound may depend.
For one observation vector, the optimal coefficient is known, and the next question is whether a single design attains it. Wu proves in \cite[Theorem~3.5]{Wu26} that every admissible scalar design has $N_\mu(R)\ge n_A(R)-O_\mu(1)$ and, for each $\varepsilon>0$, there is a design with $N_\mu(R)\le(1+\varepsilon)n_A(R)+O_\varepsilon(1)$. For an observation of rank $r$, the corresponding necessary condition is
\begin{equation}\label{eq:wu-bound}
 rN_\mu(R)\ge n_A(R)-O_\mu(1)
\end{equation}
\cite[Proposition~4.6]{Wu26}. This bound is an infinite-dimensional analogue of the dimension obstruction underlying the characterizations in \cite[Section~2]{ACMT17}. In dimension $d$, recovery from $N$ sampling times with an observation map of rank $r$ requires $rN\ge d$. Wu obtains the corresponding estimate at each spectral scale $R$: the dimension is replaced by $n_A(R)$, and observations after $R$ can compensate for only a bounded deficit in this count.

Thus, the scalar lower coefficient equals one. Our goal is to attain this lower bound up to a bounded counting error and to obtain a Riesz basis from the weighted samples.

\subsection{Interpolation and critical sampling of a single orbit}
Our first main result, \cref{thm:separated-centers},
shows how to choose points at which arbitrary
square-summable data can be interpolated uniquely and
stably by sums of translates of a fixed kernel. We consider a nonzero continuous integrable totally
positive function $\varphi$; total positivity is defined in Section~\ref{subsec:spectral-model}. The translates are centered
at an increasing sequence $(y_j)_{j\in I}$ satisfying
$y_{j+1}-y_j\ge h>0$, where $I=\{0,\ldots,N-1\}$, $I=\Nzero$,
or $I=\Z$.
For a coefficient sequence $c\in\ell^2(I)$, the corresponding
function is
\[
 F_c(u)=\sum_{j\in I}c_j\varphi(u-y_j).
\]
The interpolation problem is to choose nodes $(u_i)_{i\in I}$
so that, for every prescribed sequence $d\in\ell^2(I)$, there is
exactly one coefficient sequence $c\in\ell^2(I)$ satisfying
$F_c(u_i)=d_i$ for all $i\in I$, with a uniform bound on
$\|c\|_{\ell^2(I)}$ in terms of $\|d\|_{\ell^2(I)}$.

The theorem provides such nodes, with one node near each center.
More precisely, the nodes are increasing and satisfy
 $|u_i-y_i|\le R$, $u_{i+1}-u_i\ge\delta$,
and the collocation matrix
 $M=\bigl(\varphi(u_i-y_j)\bigr)_{i,j\in I}$
defines a boundedly invertible operator on $\ell^2(I)$.
Since the interpolation equations are exactly $Mc=d$, their
unique solution is $c=M^{-1}d$. The constants $R,\delta>0$ and
the bounds for $\|M\|$ and $\|M^{-1}\|$ depend only on $\varphi$
and $h$, not on the number or locations of the centers.
In particular, the gaps between consecutive centers need not be bounded above. A uniform coefficient estimate and an alternation argument give the finite construction; the uniform decay estimates then permit passage to infinitely many centers.

To apply the theorem to a positive single-orbit frame $(A^nf)_{n\ge0}$, we write its simple nonzero eigenvalues as $\lambda_j=e^{-s_j}$, with $s_j\downarrow0$. Their logarithmic centers $x_j=\log(1/s_j)$ are uniformly separated, and the normalized power response is a translate of
 $g(u)=\exp(u/2-e^u)$.
Applying the interpolation theorem to these centers gives a weighted real sampled orbit that is a Riesz basis. Rounding the times and adding finitely many observations gives an integer frame of finite excess. Deleting a suitable finite subset gives a weighted integer Riesz basis. Both designs satisfy
\[
 N_\mu(R)=n_A(R)+O(1).
\]
Wu's critical-count results in \cite[Theorems~3.10 and~3.11]{Wu26} assume either $x_{j+1}-x_j\to h>0$, or $x_j=jh+\eta_j$ with $h>0$ and $\eta_j$ eventually confined to an interval of length less than $h$, or $s_j/s_{j+1}\ge9$ for every $j$. Theorem~\ref{thm:critical-single} gives the critical count for every positive single-orbit frame: it uses only the separation of $(x_j)$ already implied by the frame bounds. The same designs give common sampling theorems for systems carried by a fixed separated positive spectral set.
For the same fixed pair $(A,f)$, \cref{cor:riesz-count}
gives $N_T(R)\le n_A(R)+O(1)$ for every weighted sampled
Riesz sequence. Together with the scalar lower bound,
this shows that every weighted sampled Riesz basis
satisfies $N_T(R)=n_A(R)+O(1)$.

In the positive case with $\ker A=\{0\}$, every arithmetic subsequence $(A^{Nn}f)_{n\ge0}$, for an integer $N\ge1$, remains a frame even after finitely many deletions \cite{CHPS24}. Krishtal and Miller \cite{KM26} extend the frame conclusion to $(A^{Nn+j_n}f)_{n\ge0}$ with arbitrary $j_n\in\{0,\ldots,N-1\}$. This redundancy of the unweighted orbit coexists with the weighted Riesz bases constructed here, which use the same $A$ and $f$ and have only $n_A(R)+O(1)$ observations up to time $R$. The weights make the distinction possible: $\norm{A^nf}\to0$, so no unweighted subsequence is a Riesz basis, whereas observations at the selected real or integer times are rescaled to give stable, nonredundant recovery.

The general kernel theorem also yields bases through continuous functional calculus without changing $A$ or $f$. For a continuous profile $\psi:[0,\infty)\to\R$ such that $u\mapsto e^{u/2}\psi(e^u)$ is nonzero, integrable, and totally positive, \cref{thm:ffc-critical} constructs a Riesz basis
$\bigl(\sqrt{t_k}\,\psi(t_k(I-A))f\bigr)_k$
with critical count. The matching lower bound holds for every admissible weighting within the fixed profile family. In the resolvent case, $t_j=(1-\lambda_j)^{-1}$ is an explicit choice.
Once constructed, these bases remain Riesz bases under square-summable perturbations that preserve completeness, as shown in Corollary~\ref{cor:ffc-perturbation}. For the  profiles in \eqref{eq:ffc-rational-profiles}, Proposition~\ref{prop:ffc-integer} also gives integer Riesz bases with critical count entirely within the fixed profile.

\subsection{Consequences for interpolation and related problems}
\cref{thm:separated-centers} has consequences that do not involve operator orbits.
For a fixed kernel and a fixed positive lower bound on the
separation of the centers, the selected nodes give interpolation
operators with Lebesgue constants bounded independently of the
number and locations of the centers.
The same conclusion applies to asymmetric kernels, including causal Green functions for ordinary differential equations. A further consequence concerns recovery of the weights of a location mixture from probabilities of suitably chosen disjoint intervals. We give these deductions, with their assumptions, in Section~\ref{subsub}.

In our forthcoming research, we will address several observation directions,
including optimal counts for periodic observation ranges,
geometric criteria for attainment, and questions of nonattainment
and conditioning.

\subsection{Relation to earlier work}
This paper builds on the unweighted theory of positive and sectorial orbits, which treats arithmetic sampling, redundancy, density conditions, and perturbations \cite{CHPS24,FFC26,Sectorial26,BDC26,KM26}. Related optimization problems concern the choice of generating vectors or of a common sampling support. Aguilera, Cabrelli, Negreira, and Paternostro \cite{ACNP26} study the frame index of an operator, the least number of vectors whose integer orbits form a frame. In his paper on universal sampling \cite{WuUniversal26}, Wu fixes a time set $T$ and asks whether the unweighted observations $(CA^t)_{t\in T}$ preserve observability for every positive pair $(A,C)$. Here, the operator and observation map are fixed; we choose the sampling times and one common positive weight for all observations at each time, and minimize the asymptotic number of times.

The freedom to choose weights connects our problem to scalable-frame theory \cite{KOPT13} and the scaling of dynamical frames \cite{AK17}. These weights also affect stability through the Gramian condition number $\kappa(G)=\norm G\norm{G^{-1}}$, whose optimization by diagonal scaling is studied in \cite{CC17}.

For the scalar problem, our interpolation argument builds on the classical alternating-data criteria and inverse sign patterns for totally positive matrices \cite{DFP82}, whose role in infinite-dimensional interpolation is developed in \cite{Grochenig26}. These kernels also connect sampling theory with Gabor frames \cite{Grochenig23,GRS18,GS13}. At critical density, complete interpolation is known for Gaussian and hyperbolic-secant-type shift-invariant spaces \cite{BB24,BBG22}. The sharp Kadets-type theorem in \cite[Theorem~1.2]{DGLST26} gives perturbation bounds for the sampling nodes in spaces generated by lattice translates of any continuous integrable totally positive function.

Several earlier sampling and interpolation results allow nonlattice generating centers. Gr\"ochenig and St\"ockler \cite{GS13} treat finite-type totally positive kernels with separated centers whose gaps are bounded above. Hamm and Ledford \cite{HL18} study interpolation when the generating centers form a complete interpolating sequence for a Paley--Wiener space. Ulanovskii and Zlotnikov \cite{UZAnn25} treat separated relatively dense generating sets and specified exponential-rational kernels; \cite{UZ25} concerns stability of translates along separated sets. \cref{thm:separated-centers} allows any nonzero continuous
integrable totally positive kernel and any uniformly
separated finite, one-sided infinite, or two-sided infinite
generating sequence, with no upper bound on its gaps. It selects interpolation nodes adapted to that sequence, with displacement, separation, and inverse bounds depending only on the kernel and the lower separation bound.

The node selection also relates to frame thinning and discretization. Under localization hypotheses, subframes with density arbitrarily close to critical can be selected \cite{BCHL06,BCL11}, and general discretization results for continuous frames appear in \cite{BY26,FS19}. In the present setting, the invertible collocation operator gives a bounded error between the sampling and spectral counts.

Exponential systems can also be viewed as sampled orbits. For a
measurable set $\Omega\subset\R^d$ of finite positive measure, the
operators
 $(U_tF)(\xi)=e^{2\pi i t\cdot\xi}F(\xi)$,
  $t\in\R^d$,
form a unitary representation of $\R^d$ on $L^2(\Omega)$, and
$E(\Lambda)=\{U_\lambda\mathbf{1}_\Omega:\lambda\in\Lambda\}$
is the orbit of the constant function sampled at $\Lambda$.
Fourier transformation identifies the frame inequalities for
$E(\Lambda)$ with stable sampling of functions whose Fourier
transforms are supported in $\Omega$. In this setting, criticality
is measured by Beurling density, which counts sampling parameters
in large translated regions and is compared with $|\Omega|$
\cite{BVV25,FGHKR17,Landau67}. Enstad and van Velthoven study
critical density for exponential systems \cite{EVV25} and
obstructions for smooth coherent systems \cite{EVV26}. For the
positive operator orbits considered here, the counting problem
compares the cumulative number of observation times with the
spectral count $n_A$. Exponentials also enter our proof through
Fourier transforms of translates, providing the coefficient
estimate used in the alternation argument that selects the
interpolation nodes.

The positivity assumption is essential for the scalar counting conclusion. Indeed,
Proposition~\ref{prop:normal-infinite} gives a normal operator with a
single-orbit frame whose radial spectral count has order
$R/(\log R)^2$, whereas every admissible real design requires at least
$c_\mu R/\log R$ observations up to time $R$. A different distinction
between positive and normal operator dynamics appears in the work of Krishtal
and Pfander~\cite{KP26}, who construct a normal operator with a normalized
orbit that is a frame. Their result concerns normalization of the full
integer orbit; the example here concerns the minimum number of sampling
times when arbitrary positive weights are allowed.

The remainder of the paper is organized as follows. Section~\ref{sec:prelim} collects the spectral normalization, compactness, finite completion, and rounding arguments. They are needed to connect the interpolation theorem with the construction of sampling designs that attain the critical count. Section~\ref{sec:separated} proves the interpolation theorem, explains its consequences outside dynamical sampling, and develops critical-count sampling in various settings. Finally, Section~\ref{sec:normal-obstruction} gives a normal single-orbit example for which every admissible real design has infinite counting coefficient relative to the radial spectral count.

\section{Spectral normalization, compactness, and rounding}\label{sec:prelim}
In this section, we introduce the notation and present the tools that set the stage for the main results. Spectral normalization expresses weighted operator powers
via translates of a fixed totally positive kernel,
reducing the scalar Riesz basis problem to invertibility
of a collocation operator. Uniform inverse estimates allow
the finite interpolation construction to pass to infinitely
many centers. Compact perturbation and finite completion
justify rounding the selected real sampling times to
integers. In the scalar Fredholm setting, finitely many
deletions then leave a Riesz basis. For the designs
constructed below, rounding and these finite changes
preserve the counting function up to $O(1)$.

We use $\Nzero=\{0,1,2,\ldots\}$, $\N=\{1,2,\ldots\}$, $\Z$ for the integers, and $\T=\{z\in\C:|z|=1\}$. All logarithms are natural. We write $\BB(\mathcal E,\mathcal F)$ for the space of all bounded linear maps between Hilbert spaces $\mathcal E$ and $\mathcal F$, $\sigma(A)$ for the spectrum, $\ran U$ for the range of $U$, and $I_{\mathcal E}$ for the identity on $\mathcal E$. An unsubscripted $I$ denotes the identity on the space under consideration. The notation $\ell^2(J;\mathcal E)$ means square-summable $\mathcal E$-valued families; $\ell^2(J)=\ell^2(J;\C)$. For real $v$, we write $v_+=\max\{v,0\}$.

A bounded linear map $U:\mathcal E\to\mathcal F$ is \emph{bounded below} if $\norm{Ux}\ge c\norm x$ for every $x\in\mathcal E$ and some $c>0$.
A countable family $(v_j)$ in a Hilbert space is a \emph{Bessel family} if $\sum_j|\ip{x}{v_j}|^2\le b\norm x^2$ for some $b<\infty$. It is a \emph{frame} if this sum is also at least $a\norm x^2$ for some $a>0$, and it is \emph{complete} if its closed linear span is the whole space. Its synthesis operator is $\mathcal S_vc=\sum_jc_jv_j$, the adjoint of its analysis operator. A \emph{Riesz basis} is the image of an orthonormal basis under a bounded invertible operator; equivalently, its synthesis operator is boundedly invertible. The \emph{excess} of a frame is $\dim\ker\mathcal S_v$.

\subsection{The spectral model}\label{subsec:spectral-model}
The Gramian of the integer orbit is $S=\sum_{n\ge0}A^nC^*CA^n$. Below, we refer to it as the \emph{original} Gramian. Its defining sum gives
\begin{equation}\label{eq:stein}
 C^*C=S-ASA,\qquad
 \sum_{n\ge N}A^nC^*CA^n=A^NSA^N.
\end{equation}
Therefore, for every $x\in\HH$, the lower observability bound gives
\[
 a\norm{A^Nx}^2\le\ip{SA^Nx}{A^Nx}
 =\sum_{n\ge N}\norm{CA^nx}^2\longrightarrow0.
\]
Hence $A^Nx\to0$ for every $x$. The uniform bound $\norm{A^n}^2\le b/a$ and positivity of $A$ imply $\norm A\le1$. In particular, $1$ is not an eigenvalue. For compact observation maps, these identities also force a discrete spectral structure.
\begin{theorem}\label{thm:compact-structure}
Assume that $A\in\BB(\HH)$ is positive, $C:\HH\to\YY$ is compact, and $(A,C)$ is observable. Then $A$ is a contraction, $1$ is not an eigenvalue, and $I-A$ is compact. Conversely, if $A$ is a positive contraction, $I-A$ is compact, and $\ker(I-A)=\{0\}$, then $(A,(I-A^2)^{1/2})$ is an observable system with compact observation operator and original Gramian $I$.
\end{theorem}
\begin{proof}
The contraction and stability assertions follow from \eqref{eq:stein} as described above. To prove compactness of $I-A$, we show that every spectral subspace bounded away from one is finite dimensional. We fix $0<\rho<1$ and restrict to $\HH_\rho=E_A([0,\rho])\HH$. Its compressed Gramian is
\[
 S_\rho=\sum_{n\ge0}A_\rho^n C_\rho^*C_\rho A_\rho^n,
 \qquad A_\rho=A|_{\HH_\rho},\quad C_\rho=C|_{\HH_\rho}.
\]
The series converges in operator norm because $\norm{A_\rho}\le\rho$. Every summand is compact, so $S_\rho$ is compact. Since $S_\rho\ge aI_{\HH_\rho}$, the space $\HH_\rho$ is finite dimensional. The operators $(I-A)E_A([0,\rho])$ have finite rank, and
\[
 \norm{(I-A)-(I-A)E_A([0,\rho])}
 =\norm{(I-A)(I-E_A([0,\rho]))}\le1-\rho\longrightarrow0
\]
as $\rho\uparrow1$. Thus $I-A$ is compact. Conversely, continuous functional calculus makes $(I-A^2)^{1/2}$ compact, and the spectral theorem gives $A^n\to0$ strongly. For every $N\in\Nzero$, the partial Gramians telescope:
\[
 \sum_{n=0}^{N}A^n(I-A^2)A^n=I-A^{2N+2}.
\]
The right side converges strongly to $I$. This proves observability with Gramian $I$.
\end{proof}

In infinite dimension, compact observability therefore gives the decomposition
\begin{equation}\label{eq:discrete-setting}
 \HH=K\oplus\bigoplus_{j\ge0}H_j,\qquad
 K=\ker A,\quad A|_{H_j}=e^{-s_j}I,
 \quad s_0>s_1>\cdots>0,\quad s_j\to0,
\end{equation}
where $K$ and each $H_j$ are finite dimensional. Conversely, every $C$ for which $(A,C)$ is observable in \eqref{eq:discrete-setting} is compact: $I-A$ is compact and
\[
 C^*C=(I-A)S+AS(I-A).
\]
Thus every observable pair in this discrete setting has a compact observation map. The eigenspace multiplicities may be unbounded.

We write $P_j$ for the orthogonal projection onto $H_j$. On this eigenspace, observability gives
\begin{equation}\label{eq:C-diagonal-bounds}
 a(1-e^{-2s_j})\norm{x}^2\le\norm{Cx}^2
 \le b(1-e^{-2s_j})\norm{x}^2,\qquad x\in H_j.
\end{equation}
Following \cite[Lemma~3.7]{Wu26}, we put
\begin{equation}\label{eq:Vj}
 x_j=\log(1/s_j),\qquad V_j=s_j^{-1/2}C|_{H_j},\qquad
 g(u)=\exp(u/2-e^u).
\end{equation}
Then the kernel $g$ has exponential decay
\begin{equation}\label{eq:g-decay}
 |g(u)|\le M e^{-c|u|},\qquad u\in\R,
\end{equation}
the maps $V_j$ are uniformly bounded above and below on their domains and, for $t>0$, we have
\begin{equation}\label{eq:kernel-identity}
 \sqrt t\,CA^t|_{H_j}=g(\log t-x_j)V_j.
\end{equation}
The identity \eqref{eq:kernel-identity} is the key link between
operator sampling and interpolation by translates.
For $v\in H_j$, the weighted observation $\sqrt t\,CA^tv$
depends on time through the scalar factor $g(\log t-x_j)$ only.
By \eqref{eq:g-decay}, this factor decays exponentially as
$\log t$ moves away from $x_j$, so the relevant observation
times are on the scale $t\asymp s_j^{-1}$. Moreover, according to \eqref{eq:kernel-identity}, given
  scalar observations and sampling times
$t_k=e^{u_k}$, the weighted analysis operator on
$(\ker A)^\perp$ is represented, up to a boundedly invertible
diagonal factor, by the matrix
 $M=\bigl(g(u_k-x_j)\bigr)_{k,j}$.
Thus, a boundedly invertible $M$ gives a weighted sampled orbit
that is a Riesz basis on $(\ker A)^\perp$. Moreover, when the
centers are uniformly separated, choosing one node $u_j$ within
a uniformly bounded distance of each $x_j$ gives
\[
 \#\{j:e^{u_j}\le R\}
 =\#\{j:x_j\le\log R\}+O(1)
 =n_A(R)+O(1).
\]
This explains the two requirements in the interpolation theorem
below: invertibility provides a Riesz basis, while bounded
displacement provides the critical count.

For a real function $\varphi$, its \emph{collocation matrix} at nodes $(u_i)$ and centers $(y_j)$ is $(\varphi(u_i-y_j))_{i,j}$. The function is \emph{totally positive} if $\det[\varphi(u_i-y_j)]_{i,j=1}^n\ge0$ for every $n\ge1$ and every two strictly increasing real $n$-tuples $(u_i)$ and $(y_j)$. It is \emph{strictly totally positive} if all these determinants are positive. The factorization
\begin{equation}\label{eq:strict-TP}
 g(u-y)=e^{u/2}e^{-y/2}\exp\bigl(e^u(-e^{-y})\bigr)
\end{equation}
shows that $g$ is strictly totally positive: the kernel $e^{xy}$ has this property, both changes of variable are increasing, and the remaining factors are positive \cite{Burbea74}.

\subsection{Compactness and completion}
Here, we develop compact perturbations and finite completion tools, which we later use to control the effects of rounding and restriction to one-sided sequences.

A bounded linear map $U:\mathcal E\to\mathcal F$ is \emph{Fredholm} if its range is closed and both its kernel and cokernel $\mathcal F/\ran U$ are finite dimensional. Its index is $\dim\ker U-\dim\ker U^*$. We use $\norm U_{\HS}^2=\sum_j\norm{Ue_j}^2$ for the Hilbert--Schmidt norm, where $(e_j)$ is any orthonormal basis of the domain.

For blocks $T_{ij}\in\BB(\mathcal E_j,\mathcal F_i)$, \emph{exponential off-diagonal decay} means that $\norm{T_{ij}}\le M e^{-a|i-j|}$ for some $M,a>0$. When both row and column sums below are bounded, the block matrix defines a bounded map $T:\bigoplus_j\mathcal E_j\to\bigoplus_i\mathcal F_i$ satisfying the Schur estimate
\[
 \norm T\le
 \left(\sup_i\sum_j\norm{T_{ij}}\right)^{1/2}
 \left(\sup_j\sum_i\norm{T_{ij}}\right)^{1/2}.
\]
It follows by applying Cauchy--Schwarz to each row and summing. The block Hilbert spaces need not have uniformly bounded dimension.
\begin{lemma}\label{lem:band-compact}
Let $Q=(Q_{k\ell})_{k,\ell\ge0}$ be a block matrix acting
between $\ell^2$ spaces with fixed finite-dimensional block
sizes. Suppose that $Q$ has exponential off-diagonal decay
and that
\[
 Q_{\ell+n,\ell}\longrightarrow0
 \qquad(\ell\to\infty)
\]
for every fixed $n\in\Z$, with $\ell+n\ge0$.
Then $Q$ is compact.
\end{lemma}
\begin{proof}
We first keep only finitely many diagonals. The Schur test bounds the norm of the discarded part, where $\abs{k-\ell}>M$, by a constant times $\sum_{|n|>M}e^{-c|n|}$, tending to zero. The remaining matrix is a finite sum of shifts composed with block-diagonal operators whose entries tend to zero. Each diagonal operator is compact, since truncating it to finitely many blocks gives finite-rank approximations in operator norm. Hence $Q$ is a norm limit of compact operators.
\end{proof}

We will also need off-diagonal decay for inverses. Such estimates are available for band matrices \cite{DMS84} and more general localized operators \cite{BK14}; the following elementary argument keeps the constants independent of the matrix size.

\begin{lemma}\label{lem:finite-inverse-decay}
Let $p\in[1,\infty]$, and let $B$ be an invertible matrix
indexed by a finite interval of $\Z$. Suppose that
\[
 |B_{ij}|\le M e^{-a|i-j|},
 \qquad
 \norm{B^{-1}}_{\ell^p\to\ell^p}\le K
\]
for some $M,a,K>0$.
Then there exists $b>0$, depending only on $M$, $a$, and $K$,
such that
\[
 |(B^{-1})_{ij}|\le 2Ke^{-b|i-j|}
\quad\mbox{and}\quad
 \sup_{1\le r\le\infty}
 \|B^{-1}\|_{\ell^r\to\ell^r}
 \le 2K\frac{1+e^{-b}}{1-e^{-b}},
\]
independently of the matrix size.
\end{lemma}

\begin{proof}
We use exponential conjugation, as in
\cite[Theorem~6.5]{BK14}. The idea is to weight one inverse
column by an increasing exponential and show that the
weighted column remains bounded. We denote the finite
index interval by $I\subset\Z$, fix $q\in I$, and, for
$0<b\le a/2$, set
\[
 W_b=\operatorname{diag}\bigl(e^{b|i-q|}\bigr)_{i\in I},
 \qquad B_b=W_bBW_b^{-1},
 \qquad E_b=B_b-B.
\]
Since $\bigl||i-q|-|j-q|\bigr|\le |i-j|$, we have
\begin{align*}
 |(E_b)_{ij}|
 &=|B_{ij}|
   \left|e^{b(|i-q|-|j-q|)}-1\right|
 \le Me^{-a|i-j|}\left(e^{b|i-j|}-1\right)\\
 &\le Mb|i-j|e^{-(a-b)|i-j|}
 \le Mb|i-j|e^{-a|i-j|/2},
\end{align*}
where we used $e^t-1\le te^t$ for $t\ge0$. Writing
\[
 S_a:=\sum_{m\in\Z}|m|e^{-a|m|/2}
     =\frac{2e^{-a/2}}{(1-e^{-a/2})^2},
\]
we obtain
$R_b:=\sup\limits_{i\in I}\sum\limits_{j\in I}|(E_b)_{ij}|
 \le MbS_a$
 and
 $C_b:=\sup\limits_{j\in I}\sum\limits_{i\in I}|(E_b)_{ij}|
 \le MbS_a$.
For $p\in (1,\infty)$, H\"older's inequality gives
\begin{align*}
 \|E_bx\|_{\ell^p}^p
 \le \sum_i
 \left(\sum_j|(E_b)_{ij}|\right)^{p-1}
 \sum_j|(E_b)_{ij}|\,|x_j|^p
 \le R_b^{p-1}C_b\|x\|_{\ell^p}^p.
\end{align*}
Together with the endpoint estimates
$\|E_b\|_{\ell^1\to\ell^1}\le C_b$ and
$\|E_b\|_{\ell^\infty\to\ell^\infty}\le R_b$,
this yields
$\|E_b\|_{\ell^p\to\ell^p}\le MbS_a$, $1\le p\le\infty$.

With
 $b:=\min\left\{\frac a2,\frac{1}{2KMS_a}\right\}>0$,
the assumed bound on $\|B^{-1}\|$ therefore yields
$\|B^{-1}E_b\|_{\ell^p\to\ell^p}\le\frac12$.
The factorization
 $B_b=B+E_b=B(I+B^{-1}E_b)$
and the Neumann series imply
\begin{align*}
 B_b^{-1}
 &=\sum_{n=0}^{\infty}(-B^{-1}E_b)^nB^{-1}\quad\mbox{and}
\quad\|B_b^{-1}\|_{\ell^p\to\ell^p}
 \le K\sum_{n=0}^{\infty}2^{-n}=2K.
\end{align*}
On the other hand,
$B_b^{-1}=W_bB^{-1}W_b^{-1}$ and $W_b^{-1}e_q=e_q$
give
 $\|W_bB^{-1}e_q\|_{\ell^p}
=\|B_b^{-1}e_q\|_{\ell^p}\le2K$.
For every $i\in I$, it follows that
 $e^{b|i-q|}|(B^{-1})_{iq}|
\le\|W_bB^{-1}e_q\|_{\ell^p}\le2K$.
Since $q$ was arbitrary,
\[
 |(B^{-1})_{ij}|\le2K e^{-b|i-j|},
 \qquad i,j\in I.
\]
Summing over rows or columns, we get
\[
 \max\left\{\sup_i\sum_j|(B^{-1})_{ij}|,
 \quad
 \sup_j\sum_i|(B^{-1})_{ij}|
 \right\}\le
 2K\sum_{m\in\Z}e^{-b|m|}
 =
 2K\frac{1+e^{-b}}{1-e^{-b}}.
\]
The above bound and the Schur estimate give
\[
 \sup_{1\le r\le\infty}
 \|B^{-1}\|_{\ell^r\to\ell^r}
 \le 2K\frac{1+e^{-b}}{1-e^{-b}}.
\]
All constants depend only on $M$, $a$, and $K$,
independently of $p$ and the matrix size.
\end{proof}

We will use the following compact-operator form of Peetre's lemma  \cite[Lemma~3]{Peetre61}. Suppose that $U:\mathcal E\to\mathcal F$ and
$K:\mathcal E\to\mathcal G$ are bounded operators between
Hilbert spaces, with $K$ compact. If
\begin{equation}\label{eq:peetre}
 \norm x\le M(\norm{Ux}+\norm{Kx}),
 \qquad x\in\mathcal E,
\end{equation}
then $U$ has closed range and finite-dimensional kernel. In particular, an injective operator satisfying
\eqref{eq:peetre} is bounded below. The same criterion shows that closed range and
finite-dimensional kernel are preserved under compact perturbations. Applying it also to the adjoint gives stability of Fredholmness under compact perturbations.

The preceding arguments leave at most a finite-dimensional kernel. We remove it by \emph{finite completion}: adjoining finitely many bounded observation maps so that the augmented analysis operator is bounded below. The auxiliary observations need only distinguish vectors in the remaining kernel.

\begin{lemma}\label{lem:completion}
Let $U:\HH\to\mathcal K$ be a bounded operator with closed
range and finite-dimensional kernel $E$. Let
$L_j:\HH\to\YY_j$, $j\in J$, be bounded operators into Hilbert
spaces, and suppose that
\[
 E\cap\bigcap_{j\in J}\ker L_j=\{0\}.
\]
Then there exists a finite set $F\subset J$ such that the map
\begin{equation}
    \label{eq:mapV}
    x\longmapsto\bigl(Ux,(L_jx)_{j\in F}\bigr)
\end{equation}
is bounded below. In particular, the choice $L_n=CA^n$,
$n\in\Nzero$, satisfies the hypothesis whenever $(A,C)$
is observable.
\end{lemma}
\begin{proof}
We successively choose indices $j$ so that the  maps extended by $L_j$ decrease the common kernel in $E$. In at most $\dim E$ steps their joint map $V$  given by $x\mapsto (L_jx)_{j\in F}$ is injective on $E$, hence $\norm{Vu}\ge b\norm u$ there for some $b>0$. Closed range gives $\norm{Uv}\ge c\norm v$ on $E^\perp$ for some $c >0$. For $x=u+v$ with $u\in E$ and $v\in E^\perp$,
we have
 $\norm v\le c^{-1}\norm{Ux}$  and
 $b\norm u\le\norm{Vx}+\norm V\norm v$. It follows that
 \[
\|x\|
\le
\frac1c\left(1+\frac{\|V\|}{b}\right)\|Ux\|
+\frac1b\|Vx\|
\le C\bigl(\|Ux\|^2+\|Vx\|^2\bigr)^{1/2}\quad\mbox{for some}\quad C>0.
\]
Thus, $(U,V)$ given by \eqref{eq:mapV} is bounded below.
\end{proof}

A frame whose analysis operator is Fredholm has finite excess \cite[Lemma~4.1(c)]{BCHL03}. Such a frame contains a Riesz basis obtained by deleting finitely many vectors
\cite[p.~94]{BCHL03}. Thus, in the scalar setting, finite completion of a sequence with Fredholm analysis operator followed by
a suitable finite deletion produces a Riesz basis.

\subsection{Rounding via a compact logarithm}\label{subsec:rounding}
We now turn to the passage from continuous to discrete sampling, studied for operator powers in \cite{AHP19,DMM21}. In our setting, the compactness of $\log A$ makes bounded changes of the sampling times produce a compact change in the analysis operator, as the following power-series argument shows. We keep the original indices during the perturbation, even when some times become equal, and combine the corresponding weights afterwards. This leaves the Gramian unchanged.

\begin{lemma}\label{lem:compact-time-perturbation}
Let $A$ be a bounded positive invertible operator on $\HH$,
and suppose that $B=\log A$ is compact. Let
$C:\HH\to\YY$ be bounded, and let $t_k\ge0$ and $w_k>0$ be
such that the analysis operator
\[
 Ux=(\sqrt{w_k}\,CA^{t_k}x)_k
\]
is bounded. Suppose that $\widetilde t_k\ge0$ and
\(
 \sup_k|\widetilde t_k-t_k|\le H<\infty.
\)
Then
\[
 \widetilde Ux
 =(\sqrt{w_k}\,CA^{\widetilde t_k}x)_k
\]
is bounded, and $\widetilde U-U$ is compact. If $B$ is
Hilbert--Schmidt, then $\widetilde U-U$ is Hilbert--Schmidt
as well.
\end{lemma}
\begin{proof}
 We set $\delta_k=\widetilde t_k-t_k$ and let $M_m$ multiply the $k$th output coordinate by $\delta_k^m$ so that $\norm{M_m}\le H^m$. Expanding $A^{\widetilde t_k}-A^{t_k}=A^{t_k}(e^{\delta_kB}-1) = A^{t_k}\sum_{m=1}^\infty\delta_k^m\frac{B^m}{m!}$ gives
\begin{equation}\label{eq:compact-time-series}
 \widetilde U-U=\left(\sum_{m\ge1}\frac1{m!}M_mUB^{m-1}\right)B
 \quad\mbox{and}\quad \norm{\widetilde U-U}\le\norm U(e^{H\norm B}-1),
\end{equation}
since the series in parentheses converges in operator norm. The factor $B$ makes the difference compact, and Hilbert--Schmidt when $B$ is.
\end{proof}

The above result illustrates a more general phenomenon: the essential spectral behavior of an operator orbit can determine its sampling operator up to a compact perturbation, as in the compact-localization
theorem for sectorial dynamical frames \cite{Sectorial26}.

\begin{lemma}\label{lem:logarithmic-rounding}
Let $(A,C)$ be an observable pair satisfying
\eqref{eq:discrete-setting}, with $K=\{0\}$,
$\sup_j\dim H_j<\infty$, and $x_{j+1}-x_j\ge h>0$.
Let $(u_k)_{k\ge0}$ be an increasing, bounded below,
uniformly separated sequence, and set
 $t_k=e^{u_k}$, $n_k=\lceil t_k\rceil$, and
 $Ux=(\sqrt{t_k}\,CA^{t_k}x)_k$.
If $U$ is bounded, then
\(
 U^{\rm int}x=(\sqrt{n_k}\,CA^{n_k}x)_k
\)
is bounded, and $U^{\rm int}-U$ is Hilbert--Schmidt.
No finite-dimensionality assumption on $\YY$ is needed.

If, in addition, $U$ has closed range and finite-dimensional
kernel, then the rounded family can be completed by finitely
many original integer observations. After repeated times
are combined by adding their weights, the resulting design $\nu$
satisfies
\[
 N_\nu(R)=\#\{k:t_k\le R\}+O(1).
\]
Its weights equal the positive integer sampling times
outside a finite set.
\end{lemma}
\begin{proof}
We treat the change of sampling times and the change of weights consecutively. First, the separation assumption gives $s_j\le s_0e^{-hj}$. Bounded multiplicities and \eqref{eq:C-diagonal-bounds} therefore imply
\[
 \norm{\log A}_{\HS}^2=\sum_j(\dim H_j)s_j^2<\infty
 \quad \mbox{and}\quad
 \norm C_{\HS}^2\le2b\sum_j(\dim H_j)s_j<\infty.
\]
We round the times while retaining the weights $t_k$. By \cref{lem:compact-time-perturbation}, the map $Wx=(\sqrt{t_k}\,CA^{n_k}x)_k$ is bounded and $W-U$ is Hilbert--Schmidt. It remains to replace each weight $t_k$ by $n_k$. Since $A$ is a contraction and $0\le\sqrt{n_k}-\sqrt{t_k}\le(2\sqrt{t_k})^{-1}$,
\[
 \norm{U^{\rm int}-W}_{\HS}^2
 =\sum_k(\sqrt{n_k}-\sqrt{t_k})^2\norm{CA^{n_k}}_{\HS}^2
 \le\frac{\norm C_{\HS}^2}{4}\sum_k t_k^{-1}<\infty.
\]
The last sum converges because the $u_k$ are separated and bounded below. Thus, $U^{\rm int}-U$ is Hilbert--Schmidt.

Compact perturbation and \cref{lem:completion} give finite completion. Combining repeated times preserves the Gramian. Since consecutive $t_k$ eventually differ by more than one, rounding changes counts by $O(1)$ and causes only finitely many collisions. Outside those and the completed times, the weights are $n_k$.
\end{proof}

\begin{corollary}\label{cor:real-integer-equivalence}
Let $(A,C)$ be an observable pair satisfying
\eqref{eq:discrete-setting}, with a possibly nonzero kernel
$K$. Then every admissible real design $\mu$ has an
admissible integer counterpart $\nu$ satisfying
\begin{equation}\label{eq:integer-count-comparison}
 N_\nu(R)\le N_\mu(R)+O(1),\qquad R\ge0.
\end{equation}
Consequently,
\[
 \Dens_\R(A,C)=\Dens_\Z(A,C),
\]
and the common infimum is attained in one sampling class
if and only if it is attained in the other.
\end{corollary}
\begin{proof}
We denote the projection onto $K$ by $P_K$, including $P_K=0$ when $K=\{0\}$. The operator $B=\log(A+P_K)$ is compact. We round each time upwards, retaining its weight and index. Since time zero is fixed and positive powers vanish on $K$,
\[
 A^{\lceil t\rceil}=A^t e^{(\lceil t\rceil-t)B},\qquad t\ge0.
\]
The series in \eqref{eq:compact-time-series} therefore makes the rounded analysis operator a compact perturbation of the original one. We apply finite completion and combine weights at equal integers. Rounding upwards cannot increase the count up to $R$, while completion adds only finitely many times. This proves \eqref{eq:integer-count-comparison}. Dividing by $n_A(R)\to\infty$ gives $\Dens_\Z(A,C)\le\Dens_\R(A,C)$; the reverse inequality follows because every integer design is a real design. Attainment  is transferred  via the same comparison.
\end{proof}

Thus, rounding preserves the unrestricted optimum and its attainment. For an arbitrary support, merging repeated integer times can reduce the count by an unbounded amount. With separated logarithmic sampling times, only finitely many collisions occur, giving the two-sided $O(1)$ comparison in \cref{lem:logarithmic-rounding}. In either case, however, the frame bounds and condition number may change.

\section{Interpolation and critical sampling of a single orbit}\label{sec:separated}
In this section, we first establish an interpolation theorem for totally positive kernels and then use it to construct critical sampling designs for a single orbit and for families defined through continuous functional calculus.

\subsection{Interpolation for separated centers}
We consider translates of a totally positive kernel with uniformly separated centers. The following theorem selects interpolation nodes near these centers without requiring an upper bound on the gaps.

\begin{samepage}
\begin{theorem}\label{thm:separated-centers}
\label{prop:general-tp-selection}
Let $\varphi:\R\to\R$ be a nonzero, continuous, integrable,
totally positive function, and let $h>0$. Let $I$ be $\Nzero$, $\Z$, or a
finite interval $\{0,\ldots,N-1\}$ for some $N\in\N$. Suppose that the real centers $(y_j)_{j\in I}$ satisfy
\(
 y_{j+1}-y_j\ge h\)
 whenever \(j,j+1\in I\).
Then there exist constants $R,\delta,a,b,K,\eta>0$,
depending only on $\varphi$ and $h$, and real nodes
$(u_j)_{j\in I}$ with the following properties.
\begin{enumerate}[label=\textup{(\roman*)},leftmargin=2em]
\item
The nodes satisfy
\begin{equation}\label{eq:center-displacement}
 |u_j-y_j|\le R,\qquad u_{j+1}-u_j\ge\delta
\end{equation}
whenever the indices belong to $I$.

\item
The matrix
$M=(\varphi(u_i-y_j))_{i,j\in I}$ defines a boundedly
invertible operator on $\ell^2(I)$, and
\begin{equation}\label{eq:center-frame}
 a\norm c_2^2\le\norm{Mc}_2^2\le b\norm c_2^2,\qquad
 a\norm c_2^2\le\norm{M^*c}_2^2\le b\norm c_2^2
\end{equation}
for every $c\in\ell^2(I)$.

\item
The entries of $M$ and its inverse satisfy
\[
 |M_{ij}|+|(M^{-1})_{ij}|
 \le Ke^{-\eta|i-j|},\qquad i,j\in I.
\]
\end{enumerate}
\end{theorem}
\end{samepage}
\begin{proof}
The proof has three steps. We first bound the coefficients of a finite sum of translates of a mollified kernel by its uniform norm. Secondly, we choose interpolation nodes at alternating extrema, which lets us turn this estimate into a bound for the whole inverse matrix. Finally, uniform decay allows us to pass to the original kernel and to infinitely many centers.

\emph{Uniform coefficient bounds.}
We first establish exponential decay from total positivity of order two and integrability. Total positivity implies $\varphi\ge0$, and its $2\times2$ minors give
$\varphi((u+v)/2)^2\ge\varphi(u)\varphi(v)$. Continuity then makes the positivity set an interval $J$ and $\log\varphi$ concave on $J$. If $J$ is unbounded to the right, integrability gives points $r<s$ in $J$ with $0<\varphi(s)<\varphi(r)$. For $u\ge s$, concavity bounds the slope beyond $s$ by the secant slope between $r$ and $s$, so
\[
 \varphi(u)\le\varphi(s)e^{-a_+(u-s)},\qquad
 a_+=\frac{\log\varphi(r)-\log\varphi(s)}{s-r}>0.
\]
The same argument applied to $\varphi(-u)$ gives exponential decay on an unbounded left tail. At a finite endpoint of $J$, the function vanishes beyond that endpoint. Continuity bounds the remaining compact part, and therefore
$0\le\varphi(u)\le K e^{-a|u|}$ for some $K,a>0$. In particular, $\varphi\in C_0(\R)$, the space of continuous functions vanishing at both ends of the line, and is uniformly continuous.

We use Schoenberg's representation \cite{Schoenberg51} for the second property: the Fourier transform has no real zeros. With $m=\int_\R\varphi(u)\,du>0$, the representation states that, in a vertical strip containing the imaginary axis,
\[
 \int_\R\varphi(u)e^{-zu}\,du=\frac1{\Psi(z)},\qquad
 \Psi(z)=m^{-1}e^{-\gamma z^2+d_0z}
             \prod_\nu(1+\alpha_\nu z)e^{-\alpha_\nu z},
\]
where $\gamma\ge0$, $d_0,\alpha_\nu\in\R$, and
$0<\gamma+\sum_\nu\alpha_\nu^2<\infty$. The product, which may be finite or empty, converges locally uniformly; the exponential factors cancel its linear terms. With the convention $\widehat\varphi(\xi)=\int_\R\varphi(u)e^{-iu\xi}\,du$, evaluation at $z=i\xi$ gives
\[
 |\widehat\varphi(\xi)|
 =m\exp\left(-\gamma\xi^2-\frac12\sum_\nu
                         \log(1+\alpha_\nu^2\xi^2)\right)>0.
\]
Indeed, the sum is at most $\xi^2\sum_\nu\alpha_\nu^2$. Thus the modulus of the Fourier transform has a positive minimum on every fixed compact interval.

Gaussian smoothing will give strict total positivity for the node construction. We set
\[
 G_\varepsilon(u)=(4\pi\varepsilon)^{-1/2}e^{-u^2/(4\varepsilon)},
 \qquad \varphi_\varepsilon=G_\varepsilon*\varphi,
 \qquad 0<\varepsilon\le1.
\]
Here $*$ denotes convolution, $(f*g)(u)=\int_\R f(u-v)g(v)\,dv$. We also set $\varphi_0=\varphi$. As $\varepsilon\downarrow0$, the functions $\varphi_\varepsilon$ converge uniformly to $\varphi$ and have the same $L^1$ norm. They share an exponential envelope: the original kernel has the bound already proved, while for $\varepsilon>0$ we have
\begin{equation}
 \label{eq:exp_env}
 \varphi_\varepsilon(u)\le K e^{-a|u|}
            \int_\R G_\varepsilon(v)e^{a|v|}\,dv
 \le 2K e^{a^2}e^{-a|u|}.
\end{equation}
They also share a modulus of continuity: convolution with the probability density $G_\varepsilon$ does not increase
$\|\varphi(\cdot+z)-\varphi\|_\infty$.

For $0\le\varepsilon\le1$ and a finite sum
$f=\sum_j c_j\varphi_\varepsilon(\cdot-y_j)$, we have
\[
 \widehat f(\xi)
 =e^{-\varepsilon\xi^2}\widehat\varphi(\xi)
   \sum_j c_j e^{-iy_j\xi}.
\]
We set $T=2\pi/h$. By Ingham's inequality \cite{Ingham36},
there exists $C_h>0$, depending only on $h$, such that
\[
 \int_{-T}^{T}
 \left|\sum_j c_j e^{-iy_j\xi}\right|^2\,d\xi
 \ge C_h\sum_j|c_j|^2
\]
for every finite $h$-separated family $(y_j)$ and every
choice of coefficients $(c_j)$.
Since
 $m_T:=\min_{|\xi|\le T}|\widehat\varphi(\xi)|>0$,
Plancherel's theorem gives
\begin{equation}\label{eq:tp-continuous-lower}
 \left\|\sum_jc_j\varphi_\varepsilon(\cdot-y_j)\right\|_2^2
 \ge \frac{m_T^2e^{-2T^2}}{2\pi}
       \int_{-T}^{T} \left|\sum_jc_je^{-iy_j\xi}\right|^2\,d\xi
 \ge A_{\varphi,h}\sum_j|c_j|^2,
\end{equation}
where $A_{\varphi,h}:=\frac{m_T^2e^{-2T^2}C_h}{2\pi}>0$. This constant depends only on $\varphi$ and $h$ and is uniform for $0\le\varepsilon\le1$.

The continuous Gram matrices $Q_\varepsilon$, with entries $(Q_\varepsilon)_{ij}=\int_\R\varphi_\varepsilon(u-y_i)\varphi_\varepsilon(u-y_j)\,du$, also have a uniform exponential off-diagonal envelope: integration of the common envelope bounds each entry by $C(1+|y_i-y_j|)e^{-a|y_i-y_j|}$, and $|y_i-y_j|\ge h|i-j|$. Since the lower bound \eqref{eq:tp-continuous-lower} gives
$\norm{Q_\varepsilon^{-1}}_{2\to2}\le A_{\varphi,h}^{-1}$, we may apply \cref{lem:finite-inverse-decay} with $p=2$. Its final assertion gives a uniform bound on
$\|Q_\varepsilon^{-1}\|_{\ell^\infty\to\ell^\infty}$. To pass from the $L^2$ estimate to a uniform-norm estimate, we pair $f=\sum_jc_j\varphi_\varepsilon(\cdot-y_j)$  with each translate. For every $i$,
\[
 (Q_\varepsilon c)_i
 =\int_\R f(u)\varphi_\varepsilon(u-y_i)\,du.
\]
Since $\norm{\varphi_\varepsilon}_1=\norm\varphi_1$, each pairing has modulus at most $\norm\varphi_1\norm f_\infty$. Multiplication by $Q_\varepsilon^{-1}$ gives
\begin{equation}\label{eq:tp-uniform-coefficients}
 \max_j|c_j|\le C_{\varphi,h}\norm f_\infty.
\end{equation}
All bounds so far are uniform in the number of centers and in $\varepsilon$. In the rest of the proof, constants depend only on $\varphi$ and $h$, even when only $h$ is shown as a subscript.

\emph{Choice of nodes.}
For $\varepsilon>0$, we first check that $\varphi_\varepsilon$ is strictly totally positive. The translates of $\varphi$ at distinct centers are linearly independent: after taking Fourier transforms, this follows from linear independence of distinct exponentials and the nonvanishing Fourier factor. Hence an evaluation determinant $\det[\varphi(v_i-y_j)]$ is positive on a nonempty open set of ordered $v_i$. The continuous Cauchy--Binet identity gives, for increasing $x_i,y_j$,
\[
 \det[\varphi_\varepsilon(x_i-y_j)]
 =\int_{v_1<\cdots<v_n}
 \det[G_\varepsilon(x_i-v_j)]\det[\varphi(v_i-y_j)]\,dv.
\]
The Gaussian determinant is strictly positive, while the other determinant is nonnegative and positive on that open set. Expanding the finite determinants justifies absolute convergence, so the integral is positive.

Each finite family of $n$ translates of $\varphi_\varepsilon$ is therefore a Haar system: every nonzero linear combination has at most $n-1$ distinct zeros. These functions belong to $C_0(\R)$. For $n\ge2$, we take a best uniform approximation to the last translate from the real span of the preceding $n-1$. Such an approximation exists because this span is finite dimensional. The residual is nonzero by linear independence and belongs to $C_0(\R)$. Since it tends to zero at both ends of the line, the set on which its absolute value equals its positive maximum is nonempty and compact. The alternation theorem \cite[Corollary~5.2, p.~28]{DNS80}, which applies to $C_0(I)$ for any interval $I$, including $\R$, gives increasing points $u_1<\cdots<u_n$ at which the normalized residual $f$ satisfies
$\norm f_\infty=1$ and $f(u_i)=\sigma(-1)^i$, $\sigma\in\{-1,1\}$. These extremal points will be our interpolation nodes. For $n=1$, we take a maximum point of the single translate and normalize by its value.

We write $f=\sum_jc_j\varphi_\varepsilon(\cdot-y_j)$ and form the matrix $M_\varepsilon=(\varphi_\varepsilon(u_i-y_j))$. The cofactor formula gives $(-1)^{i+j}(M_\varepsilon^{-1})_{ij}>0$. Since $M_\varepsilon c=(f(u_i))_{i=1}^n$ and $f(u_i)=\sigma(-1)^i$, the cofactor signs give
\[
 c_j=\sigma(-1)^j\sum_{i=1}^n|(M_\varepsilon^{-1})_{ji}|,\qquad 1\le j\le n.
\]
Thus applying the inverse to this alternating vector adds the absolute values in each row without cancellation. By \eqref{eq:tp-uniform-coefficients} and $\norm f_\infty=1$, we obtain
\begin{equation}\label{eq:uniform-infty-inverse}
 \norm{M_\varepsilon^{-1}}_{\infty\to\infty}
 =\max_j|c_j|\le C_{\varphi,h}.
\end{equation}
Here we have used the classical alternating-data criterion for totally positive matrices \cite{DFP82,Grochenig26}.

The common exponential envelope \eqref{eq:exp_env} gives, uniformly in $\varepsilon$, $v$, and every $h$-separated set $Y$,
\[
 \sum_{\substack{y\in Y\\|v-y|>R}}|\varphi_\varepsilon(v-y)|
 \le B_he^{-cR}.
\]
We next show that the nodes stay close to their corresponding centers. We choose $R$ so that this tail bound is less than $C_{\varphi,h}^{-1}$. If $u_k>y_k+R$, the first $k-1$ rows impose fewer conditions than there are first $k$ columns. We can therefore choose a coefficient vector of $\ell^\infty$ norm one supported on those columns and annihilating those rows. Every remaining row lies more than $R$ to the right of the chosen centers, so every remaining coordinate of the image has modulus less than $C_{\varphi,h}^{-1}$. This contradicts \eqref{eq:uniform-infty-inverse}. If $u_k<y_k-R$, we instead use the last $n-k+1$ columns and annihilate the last $n-k$ rows. The resulting coefficient vector $d$ has $\norm d_\infty=1$ and is supported on $j\ge k$. For $i\le k$ and $j\ge k$, we have $y_j-u_i\ge y_k-u_k>R$. The same tail estimate gives $\norm{M_\varepsilon d}_\infty<C_{\varphi,h}^{-1}$, again contradicting \eqref{eq:uniform-infty-inverse}. This proves the displacement bound on both sides.

It remains to bound the distances between successive nodes from below. The common envelope and translation modulus give
\begin{equation}\label{eq:tp-uniform-modulus}
 \sup_{\varepsilon,Y,v}\sum_{y\in Y}
 |\varphi_\varepsilon(v+z-y)-\varphi_\varepsilon(v-y)|\longrightarrow0
 \qquad(z\to0).
\end{equation}
For $|z|\le1$ and $L>1$, we split the sum into $|v-y|>L$ and $|v-y|\le L$. In the first region, $|v+z-y|>L-1$, so the exponential tails make the contribution uniformly small when $L$ is large. In the second region there are at most $1+2L/h$ terms, each bounded by $\norm{\varphi(\cdot+z)-\varphi}_\infty$. We first choose $L$ and then choose $|z|$ sufficiently small. This proves the uniform limit. Since consecutive alternating values differ by two, \eqref{eq:tp-uniform-coefficients} and \eqref{eq:tp-uniform-modulus} imply that successive nodes are at least a fixed distance $\delta>0$ apart. This proves \eqref{eq:center-displacement} uniformly in $\varepsilon$.

Bounded displacement gives
\begin{equation}\label{eq:finite-decay}
 |(M_\varepsilon)_{ij}|\le K_he^{-\beta_h|i-j|}.
\end{equation}
We apply \cref{lem:finite-inverse-decay} with $p=\infty$. The resulting inverse estimate and the Schur bound give both inequalities in \eqref{eq:center-frame}, uniformly for the finite matrices.

\emph{Passage to the limits.}
For a fixed finite center family, bounded displacement allows us to choose a sequence $\varepsilon_\ell\downarrow0$ along which all nodes converge. The limiting nodes retain both inequalities in \eqref{eq:center-displacement}. Uniform convergence of the kernels and convergence of the nodes imply convergence of the collocation matrices. The lower and upper bounds in \eqref{eq:center-frame} pass to the limit and show that the limit matrix is invertible. Its inverses therefore converge as well, so their exponential entry estimates also pass to the limit.

For infinitely many centers, we set
$I_n=\{0,\ldots,n-1\}$ when $I=\Nzero$ and
$I_n=\{-n,\ldots,n\}$ when $I=\Z$. The finite construction
for the original kernel $\varphi$ applies to each $I_n$
with uniform constants. Bounded displacement gives a
diagonal subsequence along which every node converges.
The limiting nodes retain both geometric bounds in
\eqref{eq:center-displacement}.

We denote the corresponding finite matrices, extended
by zero to $\ell^2(I)$, by $M_n$, and their coordinate
projections by $P_n$. For a finitely supported vector $c$,
convergence of the nodes gives coordinatewise convergence
of $M_nc$ to $Mc$ and of $M_n^*c$ to $M^*c$. The common
exponential bound \eqref{eq:finite-decay} controls the
tails of both products, so convergence holds in $\ell^2(I)$.
Once $P_nc=c$, both estimates in \eqref{eq:center-frame}
apply and pass to the limit. The Schur bound and density
extend these estimates, and the strong convergence
$M_n\to M$, to all of $\ell^2(I)$. The lower bound for $M$
gives injectivity and closed range, while the lower bound
for $M^*$ gives $\ker M^*=\{0\}$ and hence density of that
range. Thus, $M$ is onto and boundedly invertible.

The inverse decay passes to this limit as well. We extend each finite inverse by zero and denote it by $B_n$. The operators $B_n$ are uniformly bounded. For a fixed coordinate vector $e_j$, we put $x=M^{-1}e_j$. Once $P_ne_j=e_j$, the identities $B_nM_n=P_n$ and $Mx=e_j$ give
\[
 \norm{B_ne_j-P_nx}_2
 \le\norm{B_n}\,\norm{e_j-M_nx}_2\longrightarrow0.
\]
Since $P_nx\to x$, every entry of $B_n$ converges to the corresponding entry of $M^{-1}$. The uniform exponential entry bounds therefore prove (iii) for the infinite matrix.
\end{proof}

\begin{remark}\label{rem:cauchy-kernel}
For $g$ in \eqref{eq:Vj}, the continuous Gram estimate follows
directly from the Cauchy inverse formula \cite{Schechter59}.
For a finite set of centers, we put $\sigma_j=e^{-y_j}$.
The matrix indices and the product over $r\ne i$ in the
following display are restricted to this finite set. We have
\[
Q_{ij}=\int_\R g(u-y_i)g(u-y_j)\,du=\frac{\sqrt{\sigma_i\sigma_j}}{\sigma_i+\sigma_j},\qquad (Q^{-1})_{ij}=4Q_{ij}p_ip_j,\qquad p_i=\prod_{r\ne i}\frac{\sigma_i+\sigma_r}{\sigma_i-\sigma_r}.
\]
Separation yields
\[
|p_i|\le\prod_{m\ge1}\left(\frac{1+e^{-mh}}{1-e^{-mh}}\right)^2<\infty,
\qquad Q_{ij}\le e^{-h|i-j|/2},
\]
so $Q$ and $Q^{-1}$ have uniformly bounded absolute row sums. Since the Cauchy inverse formula is applied to an arbitrary finite set of centers, the finite-section lower bounds pass to the positive selfadjoint infinite matrix, making it boundedly invertible as well. The bounds for $Q^{-1}$ and $\norm g_1$ give \eqref{eq:tp-uniform-coefficients} directly. Since $g$ is strictly totally positive, Gaussian regularization is unnecessary; the proof proceeds with the node construction and limit argument above.
\end{remark}

For the center set $Y=\{y_j:j\in I\}$ in \cref{thm:separated-centers}, we obtain complete interpolation in
\[
 V_Y^2(\varphi)=\left\{\sum_{j\in I}c_j\varphi(\cdot-y_j):c\in\ell^2(I)\right\},
\]
meaning that evaluation $F\mapsto(F(u_i))_{i\in I}$ is a bounded bijection from this space onto $\ell^2(I)$. The series converge in $L^2(\R)$ and uniformly on $\R$, so their sums are continuous; the $L^2$ and coefficient norms are equivalent by  \eqref{eq:tp-continuous-lower} and the exponential envelope. Complete interpolation holds at the selected nodes; for other bounded displacements, invertibility must be checked separately. The theorem applies, for example, to Gaussian kernels, $e^{-|u|}$, and $g$ in \eqref{eq:Vj}.

\subsubsection{Other consequences of Theorem~\ref{thm:separated-centers}}\label{subsub}
The finite-dimensional conclusion of \cref{thm:separated-centers} also gives consequences outside dynamical sampling. We describe three of them. Throughout, the constants may depend on the fixed kernel $\varphi$ and the lower separation $h$, but not on the number or locations of the centers.

First, we consider interpolation by translates of $\varphi$.
Uniform bounds on the associated Lebesgue constants control
both the amplification of errors in the data and the loss
relative to best approximation; see, for example, \cite{HNW10}.
For every finite $h$-separated center set, our theorem provides
nearby interpolation nodes for which these bounds depend only
on $\varphi$ and $h$.
We work in the space
\[
 V_Y=\operatorname{span}\{\varphi(\cdot-y_k):0\le k<N\}.
\]
For the selected nodes, we put $M_{ik}=\varphi(u_i-y_k)$ and define
\[
 L_j(x)=\sum_{k=0}^{N-1}(M^{-1})_{kj}\varphi(x-y_k),\qquad
 P_Yf(x)=\sum_{j=0}^{N-1}f(u_j)L_j(x).
\]
Then $L_j(u_i)=\delta_{ij}$, and $P_Yv=v$ for every $v\in V_Y$. The inverse estimate gives
$\sup_k\sum_j|(M^{-1})_{kj}|\le C_{\varphi,h}$. Separation and the exponential envelope give
$\sup_x\sum_k|\varphi(x-y_k)|\le C'_{\varphi,h}$. Consequently,
\[
 \sup_{x\in\R}\sum_j|L_j(x)|\le\Lambda_{\varphi,h}.
\]
Thus the Lebesgue constants are uniformly bounded, and every bounded continuous $f$ satisfies
\[
 \norm{f-P_Yf}_\infty
 \le(1+\Lambda_{\varphi,h})\inf_{v\in V_Y}\norm{f-v}_\infty.
\]
Moreover, if $B_L$ is obtained from $M^{-1}$ by setting its entries with $|i-j|>L$ equal to zero, the Schur estimate gives
$\norm{M^{-1}-B_L}_{2\to2}\le Ce^{-cL}$.
Thus, truncation to $L$ diagonals on either side of the main
diagonal approximates $M^{-1}$ with an exponentially small
operator-norm error, uniformly in the number of centers.

Second, we consider the causal kernel $\varphi(x)=x_+e^{-x}$. It is the convolution of two copies of $e^{-x}\mathbf1_{[0,\infty)}(x)$, and is therefore continuous, integrable, and totally positive \cite{Schoenberg51}. In the sense of distributions,
\[
 (D+1)^2\varphi=\delta_0,\qquad D=\frac{d}{dx}.
\]
For prescribed separated source locations $y_j$, the function
$F=\sum_jc_j\varphi(\cdot-y_j)$ is the causal solution of
\[
 (D+1)^2F=\sum_jc_j\delta_{y_j}.
\]
The theorem supplies one observation point per source at which the source strengths $c_j$ are recovered with uniform $\ell^2$ stability. Evaluation at the source locations themselves would give the strictly lower triangular matrix $(\varphi(y_i-y_j))$, which is singular. The same interpretation applies to
$\varphi_m(x)=x_+^{m-1}e^{-x}/(m-1)!$, $m\ge2$, the causal Green functions of $(D+1)^m$. 

This example is related to the spline-recovery problem
discussed in \cite{KP25}. Indeed, multiplication by $e^x$
transforms
$F_m=\sum_jc_j\varphi_m(\,\cdot-y_j)$ into a polynomial spline whose $m$th
distributional derivative is
$\sum_j c_je^{y_j}\delta_{y_j}$. Here, the source locations are prescribed, and the conclusion concerns the uniformly stable recovery of their amplitudes from one point
observation per source.

Finally, we consider estimation of mixture weights when the
component densities are known; see, for example, \cite{KCSA20}.
For translates of a fixed probability density, the theorem
provides observation intervals whose bin probabilities determine
the weights with stability bounds independent of the number
of components, provided the centers have a fixed positive
separation.
We normalize $\int_\R\varphi=1$ and consider the probability density
\[
 p_w(x)=\sum_{j=0}^{N-1}w_j\varphi(x-y_j),\qquad
 w_j\ge0,\quad \sum_jw_j=1,
\]
with known locations and unknown weights. If $0<\rho<\delta/2$, the intervals
$I_i=[u_i-\rho,u_i+\rho]$ are disjoint. We define
$Q_{ij}=\int_{I_i}\varphi(x-y_j)\,dx$. Uniform continuity, exponential decay, and separation of both the nodes and centers imply
\[
 \left\|(2\rho)^{-1}Q-M\right\|_{2\to2}\longrightarrow0
 \qquad(\rho\downarrow0)
\]
uniformly in $N$ and the center configuration. Indeed, for the absolute row and column sums we first bound the terms with $|u_i-y_j|>L$ by the exponential tails. There are uniformly boundedly many remaining terms in each row and column, and their differences tend uniformly to zero by continuity. The Schur estimate proves the assertion.

We may therefore fix $\rho>0$, depending only on $\varphi$ and $h$, so that $Q$ is invertible with a uniformly bounded inverse. The bin probabilities satisfy $q=Qw$. For $m$ independent observations, the empirical bin frequencies $\widehat q$ satisfy
$\mathbb E\norm{\widehat q-q}_2^2=m^{-1}\sum_iq_i(1-q_i)\le m^{-1}$. If $\widehat w$ is the Euclidean projection of $Q^{-1}\widehat q$ onto the probability simplex, its nonexpansiveness gives
\[
 \mathbb E\norm{\widehat w-w}_2^2\le C_{\varphi,h}/m.
\]
This is an $\ell^2$ estimate for known-location mixture weights. It does not assert recovery of unknown locations or a dimension-independent bound in total variation.

\subsection{Critical sampling of a single orbit}
For power observations, the response $e^{-s_jt}$ gives the kernel $g$ in \eqref{eq:Vj} through \eqref{eq:kernel-identity}. We first verify that separated scalar spectra give observable integer orbits, in other words, Carleson frames, and then construct critical real and integer designs. Here $(e_j)$ denotes an orthonormal eigenbasis and $f_j=\ip f{e_j}$ is the corresponding coordinate of $f$.

\begin{lemma}\label{lem:canonical-separated}
Let $(e_j)_{j\ge0}$ be an orthonormal basis of $\HH$, and
let $(s_j)_{j\ge0}$ be a positive sequence such that
\(s_j\downarrow0\) and
 \(\log(s_j/s_{j+1})\ge h>0\).
Define
\[
 Ae_j=e^{-s_j}e_j,\qquad
 f=\sum_{j\ge0}b_j\sqrt{s_j}\,e_j,
\]
where
\(0<\inf_j|b_j|\le\sup_j|b_j|<\infty\).
Then $(A^nf)_{n\ge0}$ is a frame.
\end{lemma}
\begin{proof}
We put $\lambda_j=e^{-s_j}$. Since $s_{j+1}/s_j\le e^{-h}$ and
$(1-e^{-s})/s\to1$ as $s\downarrow0$, we have
\[
 \limsup_{j\to\infty}
 \frac{1-\lambda_{j+1}}{1-\lambda_j}
 \le e^{-h}<1.
\]
Each of these ratios is strictly less than one, so their
supremum is also less than one. Hence $(\lambda_j)$ satisfies
the Carleson condition by \cite[Proposition~1.2(i)]{CHPS24}.
Moreover,
\[
 f=\sum_{j\ge0}m_j\sqrt{1-\lambda_j^2}\,e_j,
 \qquad
 m_j=b_j\sqrt{\frac{s_j}{1-e^{-2s_j}}},
\]
where $0<\inf_j|m_j|\le\sup_j|m_j|<\infty$.
The conclusion now follows from \cite[Theorem~1.1]{CHPS24}.
\end{proof}

We can now select critically sparse sampling times for a single orbit. The proof derives the required spectral separation directly from the frame bounds. The complete characterization of single normal-orbit frames is given in \cite{ACCMP17,ACMT17}. The integer construction has only finite excess, so finitely many deletions leave a Riesz basis \cite{BCHL03}.

\begin{theorem}\label{thm:critical-single}
Let $A\ge0$ be a bounded operator on an infinite-dimensional
Hilbert space $\HH$, and let $f\in\HH$. Suppose that
$(A^nf)_{n\ge0}$ is a frame. Then there exist a real design
and an integer design whose weighted sampled orbits are
Riesz bases. Each of these designs $\mu$ satisfies
\[
 N_\mu(R)=n_A(R)+O(1).
\]
For the integer design, the weights satisfy $w_t=t$ at all
but finitely many positive sampling times.
\end{theorem}
\begin{proof}

We apply \cref{thm:compact-structure} to $Cx=\ip{x}{f}$. Observability makes $C$ injective on each eigenspace. Hence $\dim\ker A\le1$, and the nonzero eigenvalues are simple and accumulate only at one. We enumerate them as $\lambda_j=e^{-s_j}\uparrow1$ and choose an orthonormal eigenbasis $(e_j)$ of $(\ker A)^\perp$. We write $f_j=\ip f{e_j}$. Equation~\eqref{eq:C-diagonal-bounds} gives $f_j=b_j\sqrt{s_j}$, where $|b_j|$ is bounded above and away from zero.

 Next, we derive separation from the original frame bounds $a,b$. For $i\ne j$, we compress the frame operator to $\operatorname{span}\{e_i,e_j\}$ and denote its matrix by $G$. We have $aI\le G\le bI$. With $D=\diag(G_{11},G_{22})$, the normalized matrix is
\[
 D^{-1/2}GD^{-1/2}=\begin{pmatrix}1&z\\\overline z&1\end{pmatrix}
 \ge(a/b)I.
\]
Its eigenvalues are $1\pm|z|$, so $|z|\le1-a/b$. Direct calculation gives
\[
 |z|=\frac{\sqrt{(1-\lambda_i^2)(1-\lambda_j^2)}}{1-\lambda_i\lambda_j}
 =\sqrt{1-\rho(\lambda_i,\lambda_j)^2},\qquad
 \rho(\lambda_i,\lambda_j)=\frac{|\lambda_i-\lambda_j|}{1-\lambda_i\lambda_j}.
\]
Since $|z|>0$, we have $a<b$. Thus
$\rho(\lambda_i,\lambda_j)\ge\delta:=\sqrt{1-(1-a/b)^2}\in(0,1)$. For $0<s<u$, the monotonicity of $\sinh v/v$ on $(0,\infty)$ gives
\[
 \rho(e^{-s},e^{-u})
 =\frac{\sinh((u-s)/2)}{\sinh((u+s)/2)}
 \le\frac{u-s}{u+s}.
\]
Taking $s=s_{j+1}$ and $u=s_j$, we obtain
\[
 x_{j+1}-x_j=\log(s_j/s_{j+1})\ge
 h:=\log\frac{1+\delta}{1-\delta}>0,\qquad x_j=\log(1/s_j).
\]

We apply \cref{thm:separated-centers} with $\varphi=g$ to these centers, where $g$ is defined in \eqref{eq:Vj}. We put $t_j=e^{u_j}$ and assign weight $t_j$. In the eigenbasis above, the analysis operator on $(\ker A)^\perp$ is
$U_+=M\diag(\overline{b_j})$, where $M=(g(u_i-x_j))_{i,j\ge0}$. Both factors are boundedly invertible. Hence, the weighted vectors form a Riesz basis for this subspace.

If $\ker A\ne \{0\}$, we choose a unit vector $e_*$ spanning it. Observability gives $\gamma=\ip f{e_*}\ne0$. We add time zero with weight one. With this observation listed first, the augmented analysis operator on $\C e_*\oplus(\ker A)^\perp$ has the form
\[
 \begin{pmatrix}\overline\gamma&L\\0&U_+\end{pmatrix},
 \qquad Lx=\ip x f\quad(x\in(\ker A)^\perp).
\]
It is boundedly invertible, so the augmented family is a Riesz basis for the whole space. Bounded displacement and separation give
\[
 \#\{j:u_j\le\log R\}=\#\{j:x_j\le\log R\}+O(1)=n_A(R)+O(1).
\]
The possible time-zero observation changes this count by at most one.

For the integer design, we apply \cref{lem:logarithmic-rounding} on $(\ker A)^\perp$ and retain time zero when needed. Before coincident times are merged, the rounded analysis operator is a Hilbert--Schmidt perturbation of the invertible real analysis operator; hence it is Fredholm. Merging the finitely many coincident times preserves Fredholmness. Finite completion makes the analysis operator bounded below while preserving its finite-dimensional cokernel. The resulting weighted family is, therefore, a frame of finite excess. The finite-deletion argument following \cref{lem:completion} leaves a Riesz basis. Rounding, completion, and deletion change the count by only $O(1)$ and affect only finitely many of the weights $w_t=t$. Finally, Wu's lower bound \eqref{eq:wu-bound}, with rank one, proves optimality.
\end{proof}

\begin{corollary}\label{cor:riesz-count}
Let $A$ and $f$ satisfy the hypotheses of
\cref{thm:critical-single}. If a weighted sampled orbit
$(\sqrt{w_t}A^tf)_{t\in T}$ is a Riesz sequence, then
 $N_T(R)\le n_A(R)+O(1)$.
If it is a Riesz basis, then
 $N_T(R)=n_A(R)+O(1)$.
 The constants in the bounded errors may depend on $A$,
$f$, and the Riesz bounds of the sampled family, but
not on $R$.
\end{corollary}

\begin{proof}
As in the proof of \cref{thm:critical-single}, we enumerate
the nonzero eigenvalues of $A$ as $e^{-s_j}$, $j\ge0$,
where $s_j\downarrow0$, and choose corresponding
orthonormal eigenvectors $e_j$. With $f_j=\ip f{e_j}$,
there exist constants $c,C>0$ and $q\in(0,1)$ such that
 $c s_j\le |f_j|^2\le C s_j$,
$s_{j+1}\le q s_j$.
Moreover, $d:=\dim\ker A\le1$.

We write $v_t=\sqrt{w_t}A^tf$ and denote its Riesz
sequence bounds by $\alpha,\beta>0$. In particular,
$(v_t)_{t\in T}$ is a Bessel family with bound $\beta$.
For $R>0$, we put
\[
 m=n_A(R)=\#\{j:s_j\ge R^{-1}\},
 \qquad W(R)=\sum_{\substack{t\in T\\t\le R}}w_t.
\]
Since $Rs_m<1$, the Bessel inequality applied to $e_m$
gives
\[
 \beta
 \ge |f_m|^2
       \sum_{\substack{t\in T\\t\le R}}w_t e^{-2ts_m}
 \ge c s_m e^{-2}W(R).
\]
Consequently,
 $W(R)\le \frac{\beta e^2}{c s_m}$.

For an integer $L\ge1$, we denote by $P$ the orthogonal
projection onto
\[
 \ker A\oplus
 \operatorname{span}\{e_0,\ldots,e_{m+L-1}\}.
\]
The preceding estimate and geometric decay yield
\[
 \begin{aligned}
 \sum_{\substack{t\in T\\t\le R}}\|(I-P)v_t\|^2
 &=\sum_{\substack{t\in T\\t\le R}}
       w_t\sum_{j\ge m+L}|f_j|^2e^{-2ts_j}\le W(R)\,C\sum_{j\ge m+L}s_j\le \frac{\beta e^2C}{c(1-q)}q^L.
 \end{aligned}
\]
We fix $L$, independently of $R$, so that the last
expression is less than $\alpha$.

If $N_T(R)>m+L+d=\rank P$, there are coefficients
$(a_t)_{t\in T\cap[0,R]}$ satisfying
 $\sum_{t\le R}|a_t|^2=1$
 and $P\sum_{t\le R}a_tv_t=0$.
The lower Riesz bound and Cauchy--Schwarz then imply
\[
 \alpha
 \le \left\|\sum_{t\le R}a_tv_t\right\|^2
 =\left\|\sum_{t\le R}a_t(I-P)v_t\right\|^2
 \le\sum_{t\le R}\|(I-P)v_t\|^2
 <\alpha,
\]
a contradiction. Thus,
 $N_T(R)\le n_A(R)+L+d$.

If the sampled orbit is a Riesz basis, its design is
admissible. The rank-one case of \eqref{eq:wu-bound}
gives $N_T(R)\ge n_A(R)-O(1)$, completing the proof.
\end{proof}

The same counting identity holds for every weighted sampled
frame of finite excess for the fixed pair $(A,f)$. Indeed,
such a frame contains a Riesz basis obtained by deleting
finitely many vectors \cite[p.~94]{BCHL03}, and these
deletions change the counting function by only $O(1)$.

A design is \emph{common} to a class of pairs if the same support and weights are admissible for every pair in the class; the frame bounds may depend on the pair. The scalar construction gives common real and integer designs for systems whose spectral measures are carried by a prescribed set $\Lambda$, meaning $E_A(\R\setminus\Lambda)=0$. The eigenspaces and output spaces may be infinite-dimensional.

\begin{corollary}\label{cor:common-critical}
Let $\Lambda=\{e^{-s_j}:j\ge0\}$, where $s_j\downarrow0$ and
\(\inf_j\log(s_j/s_{j+1})>0\),
and set
\(m_\Lambda(R)=\#\{j:s_j\ge R^{-1}\}\).
Then there exist real and integer designs admissible for
every observable positive pair $(A,C)$ whose spectral
measure is carried by $\Lambda$. Each of these designs satisfies
 $N_\mu(R)=m_\Lambda(R)+O(1)$.
The observation space and the eigenspaces may be
infinite-dimensional. Every design common to this class
satisfies
\(N_\mu(R)\ge m_\Lambda(R)-O_\mu(1)\).
\end{corollary}
\begin{proof}
The canonical pair on $\ell^2(\Nzero)$ given by
\[
 A_0e_j=e^{-s_j}e_j,\qquad
 f^{(0)}=\sum_{j\ge0}\sqrt{s_j}\,e_j,\qquad
 C_0x=\ip{x}{f^{(0)}}
\]
is observable by \cref{lem:canonical-separated}.
Theorem~\ref{thm:critical-single} provides real and integer
designs for this pair with counting function
$m_\Lambda(R)+O(1)$.
By \cite[Theorem~4.9]{Wu26}, both designs are admissible
for every pair in the stated class.

Conversely, every common design is admissible for the
canonical pair. The rank-one case of \eqref{eq:wu-bound}
therefore gives
$N_\mu(R)\ge m_\Lambda(R)-O_\mu(1)$.
\end{proof}

\begin{example}\label{ex:critical-support-instability}\label{rem:evv-deformation}
A critical support need not remain admissible under small changes of the operator and generating vector. The spectral count can change even though each perturbed integer orbit is still a frame.

For a fixed $h>0$ and any $\tau\in(0,1]$, we set
\[
 s_j(\tau)=e^{-\tau hj},\qquad
 A_\tau e_j=e^{-s_j(\tau)}e_j,\qquad
 f_\tau=\sum_{j\ge0}\sqrt{s_j(\tau)}\,e_j.
\]
Each integer orbit is a frame by \cref{lem:canonical-separated}, and
 $\|A_\tau-A_1\|+\|f_\tau-f_1\|\longrightarrow0$ $(\tau\uparrow1)$.
For the first limit, we use convergence on finitely many coordinates and the uniform estimate $|e^{-s_j(\tau)}-1|\le e^{-\tau_0hj}$ on the tail when $\tau\ge\tau_0>0$. For the second, we use the summable bound $4e^{-\tau_0hj}$ for the squared coordinate differences. We choose a support $T$ with critical-count attainment for $(A_1,f_1)$ and write $N_T(R)=\#(T\cap[0,R])$. Then
\[
 N_T(R)=h^{-1}\log R+O(1),\qquad
 n_{A_\tau}(R)=(\tau h)^{-1}\log R+O(1).
\]
The scalar lower bound rules out every admissible weighting on $T$ for every $\tau<1$. Although the operator and sensor converge in norm, the logarithmic centers move by $-(1-\tau)hj$, an unbounded displacement. Norm convergence of the operator and sensor therefore does not give the small  perturbation of the analysis operator required for frame stability.
\end{example}

\subsection{Critical sampling through functional calculus}\label{subsec:ffc}
The same interpolation argument applies when powers of $A$ are replaced by certain other continuous functions. Keeping $A$ and the generating vector $f$ fixed, we seek bases with the same critical count. Here $t$ serves as a scale parameter for the chosen family of functions.

We assume that $A$ and $f$ satisfy \cref{thm:critical-single} and enumerate \emph{all} eigenvalues as $0\le\lambda_0<\lambda_1<\cdots<1$, including the possible zero eigenvalue. We put
\[
 B=I-A,\qquad r_j=1-\lambda_j,\qquad
 y_j=\log(1/r_j),\qquad f_j=d_j\sqrt{r_j}.
\]
The original frame bounds imply $a\le|d_j|^2\le2b$. We write $d_-=\inf_j|d_j|>0$ and $d_+=\sup_j|d_j|<\infty$. The normalized two-dimensional Gramian calculation in the proof of \cref{thm:critical-single} gives a pseudohyperbolic lower bound $\delta>0$ for every pair of distinct eigenvalues, including a possible zero eigenvalue. For $r_i>r_j$, we have
\[
 \delta\le\rho(\lambda_i,\lambda_j)
 =\frac{r_i-r_j}{r_i+r_j-r_ir_j}
 \le1-\frac{r_j}{r_i},
\]
because $r_i\le1$. Thus $r_j/r_i\le1-\delta$, and the centers satisfy $y_{j+1}-y_j\ge h_1:=-\log(1-\delta)>0$. We shall use the comparison
\begin{equation}\label{eq:ffc-count-comparison}
 \#\{j:r_j\ge R^{-1}\}=n_A(R)+O(1).
\end{equation}
To see this, we omit a possible zero eigenvalue. On the remaining spectrum, $r_j\le-\log\lambda_j\le K_A r_j$ for a fixed $K_A<\infty$. Separation bounds the number of $y_j$ in every interval of length $\log K_A$, so changing the threshold by this factor changes the count by $O(1)$.

For $\psi\in C([0,\infty);\R)$, we define
\begin{equation}\label{eq:ffc-profile}
 \varphi(u)=e^{u/2}\psi(e^u),\qquad
 F_t(z)=\sqrt t\,\psi\bigl(t(1-z)\bigr),\quad t>0.
\end{equation}
For each fixed $t$, $F_t$ belongs to $C([0,1])$, the space of continuous scalar functions with the uniform norm, and $F_t(A)$ is a bounded operator. Spectral calculus gives the exact identity
\begin{equation}\label{eq:ffc-kernel-identity}
 F_t(A)f=\sqrt t\,\psi(tB)f
       =\sum_j d_j\varphi(\log t-y_j)e_j.
\end{equation}
For a design $\mu=\sum_{t\in T}w_t\delta_t$, we call the family
$(\sqrt{w_t}\,\psi(tB)f)_{t\in T}$ admissible when it is a frame. Thus the profile $\psi$ specifies the family whose sampling count we minimize. At $t=0$ the vector is $\sqrt{w_0}\,\psi(0)f$.

\begin{theorem}\label{thm:ffc-critical}
Let $A$ and $f$ satisfy the hypotheses of
\cref{thm:critical-single}. Let $(\lambda_k)_{k\ge0}$
enumerate all eigenvalues of $A$ in increasing order,
including a possible zero eigenvalue, and let
$\psi\in C([0,\infty);\R)$. Suppose that the function
$\varphi$ in \eqref{eq:ffc-profile} is nonzero, integrable,
and totally positive.

Then there exists an increasing sequence of positive
parameters $(t_k)_{k\ge0}$ with uniformly separated
logarithms and
$\sup_{k\ge0}\bigl|\log\bigl(t_k(1-\lambda_k)\bigr)\bigr|
 <\infty$
such that, for $T=\{t_k:k\ge0\}$,
\begin{equation}\label{eq:ffc-riesz}
 \bigl(\sqrt{t_k}\,\psi(t_k(I-A))f\bigr)_{k\ge0}
 \quad\text{is a Riesz basis and}\quad
 N_T(R)=n_A(R)+O(1).
\end{equation}
Moreover, every admissible weighted design $\mu$ for the
profile $\psi$ satisfies
\begin{equation}\label{eq:ffc-lower-count}
 N_\mu(R)\ge n_A(R)-O_\mu(1).
\end{equation}
Thus, the optimal real-parameter counting coefficient
relative to $n_A$ is one, and it is attained with bounded
counting error. The lower bound remains valid when
finitely many additional vectors are allowed, with their
number absorbed into the constant.
\end{theorem}

\begin{proof}
We apply \cref{thm:separated-centers} to $\varphi$ and  \(y_j=-\log(1-\lambda_j)\), obtaining nodes $u_k$. We put $t_k=e^{u_k}$. By \eqref{eq:ffc-kernel-identity}, the analysis operator is
\[
 \bigl(\varphi(u_k-y_j)\bigr)_{k,j}\diag(\overline{d_j}),
\]
which is invertible. Bounded logarithmic displacement, separation, and \eqref{eq:ffc-count-comparison} give the count in \eqref{eq:ffc-riesz}.

For the lower count, we show that observations after $R$ carry too little energy on all but a fixed number of the first $n_A(R)$ spectral modes. Observations up to $R$ must then distinguish those modes, forcing the count. The needed tail estimate comes from total positivity of order two and exponential decay on the right. We choose $v_0$ with $\varphi(v_0)>0$. The minor with row points $v_0+d,v+d$ and column points $0,d$ gives
\[
 \varphi(v+d)\varphi(v_0)
 \le\varphi(v)\varphi(v_0+d),\qquad v\ge v_0,\ d\ge0.
\]
Together with the exponential envelope used in the proof of \cref{thm:separated-centers}, this yields constants $L,\gamma>0$ such that
\begin{equation}\label{eq:ffc-one-sided-ratio}
 0\le\varphi(v+d)\le L e^{-\gamma d}\varphi(v),
 \qquad v\ge v_0,\ d\ge0.
\end{equation}
The inequality remains valid when $\varphi(v)$ is zero.

We denote the sampled frame bounds by $\alpha,\beta$ and set
$m=\#\{j:y_j\le\log R-v_0\}$. For $m>K$, $t>R$, and $0\le j<m-K$, we apply \eqref{eq:ffc-one-sided-ratio} with $v=\log t-y_{m-1}$ and $d=y_{m-1}-y_j$. Equation \eqref{eq:ffc-kernel-identity} gives
\[
 |f_j\psi(tr_j)|^2
 \le \left(\frac{d_+}{d_-}\right)^2 L^2
 e^{-2\gamma(y_{m-1}-y_j)}|f_{m-1}\psi(tr_{m-1})|^2,
\]
where $d_-\le|d_j|\le d_+$. Testing the upper frame bound on $e_{m-1}$ and summing over $j<m-K$ bounds the trace of the late-time Gramian on
$E=\operatorname{span}\{e_0,\ldots,e_{m-K-1}\}$ by
\begin{equation}
    \label{eq:trace_bound}
 \beta\left(\frac{d_+}{d_-}\right)^2 L^2
 \sum_{\ell\ge K}e^{-2\gamma h_1\ell}.
\end{equation}
 Comparison with $e_{m-1}$ controls the entire sum
by one Bessel bound, while separation makes the displayed
tail uniformly small for large $K$. Since the compressed
Gramian is positive, its trace bounds the late-time
energy of every unit vector in $E$.
We choose $K$, independently of $R$, so that the trace bound \eqref{eq:trace_bound} is less than $\alpha$. If $N_\mu(R)<m-K=\dim E$, the observations at $t\le R$ leave a unit vector of $E$ undetected. Its late-time energy is bounded by the trace estimate, hence is also less than $\alpha$, contradicting the lower frame bound. By separation and \eqref{eq:ffc-count-comparison}, $m=n_A(R)+O(1)$, proving \eqref{eq:ffc-lower-count}. With $q$ additional vectors, we impose their $q$ linear conditions as well; the same argument gives $N_\mu(R)\ge m-K-q$.
\end{proof}

\subsubsection{The value at the spectral accumulation point.}
Starting from a totally positive kernel $\varphi$, the corresponding profile is
$\psi(s)=s^{-1/2}\varphi(\log s)$ for $s>0$. Continuity at zero is the remaining condition for bounded continuous functional calculus. It depends on the decay of $\varphi$ at $-\infty$: for example, $\varphi(u)=e^{-\varepsilon|u|}$ with $0<\varepsilon<1/2$ gives $s^{\varepsilon-1/2}$ on $(0,1]$, which is unbounded at zero. Dilating the kernel increases its decay rate and supplies the required continuity.

\begin{corollary}\label{cor:ffc-dilation}
Let $A$ and $f$ satisfy the hypotheses of
\cref{thm:critical-single}, and let $\varphi:\R\to\R$ be
a nonzero, continuous, integrable, totally positive
function. Suppose that
\(
 |\varphi(u)|\le Ke^{-a|u|}\), \(u\in\R\),
for some $K, a>0$. For $\tau>1/(2a)$, define
\[
 \psi_\tau(s)=s^{-1/2}\varphi(\tau\log s)\quad(s>0),
 \qquad \psi_\tau(0)=0.
\]
Then $\psi_\tau$ satisfies the profile assumptions of
\cref{thm:ffc-critical}. Consequently, every such kernel,
after a fixed dilation if needed, yields a Riesz basis
through continuous functional calculus for the same
$A$ and $f$, with sampling count $n_A(R)+O(1)$.
\end{corollary}
\begin{proof}
Dilation preserves total positivity and integrability. For $0<s\le1$, we have
$|\psi_\tau(s)|\le K s^{a\tau-1/2}\to0$, while
$e^{u/2}\psi_\tau(e^u)=\varphi(\tau u)$. We apply \cref{thm:ffc-critical}. The exponential envelope follows from log-concavity and integrability, as shown in the proof of \cref{thm:separated-centers}.
\end{proof}

\subsubsection{Resolvents and other profiles.}
Sampling formulas associated with operator resolvents have
been studied through resolvent-type kernels and their
reproducing kernel spaces; see \cite{FGH13}.
In the present setting, the eigenvalues determine an explicit
family of resolvent vectors forming a Riesz basis with optimal
parameter count.

\begin{proposition}\label{prop:ffc-resolvent}
Let $A$ and $f$ satisfy the hypotheses of
\cref{thm:critical-single}. Let $(\lambda_j)_{j\ge0}$
enumerate all eigenvalues of $A$ in increasing order,
including a possible zero eigenvalue, and set
$r_j=1-\lambda_j$. Then
\begin{equation}\label{eq:explicit-resolvent-basis}
 \left(
 \sqrt{r_j}\bigl((1+r_j)I-A\bigr)^{-1}f
 \right)_{j\ge0}
\end{equation}
is a Riesz basis.
Equivalently, this family corresponds to $\psi(s)=(1+s)^{-1}$ and $t_j=r_j^{-1}$ in \eqref{eq:ffc-riesz}. Its parameter count is $\#\{j:r_j\ge R^{-1}\}=n_A(R)+O(1)$ and is optimal among the weighted designs for this profile.
\end{proposition}
\begin{proof}
The profile kernel is
$\varphi(u)=e^{u/2}/(1+e^u)=1/(2\cosh(u/2))$. The Cauchy determinant, after positive row and column scalings, gives its total positivity; it is continuous and integrable. With $t_i=r_i^{-1}$, the analysis matrix of \eqref{eq:explicit-resolvent-basis} is
\[
 Q\diag(\overline{d_j}),\qquad
 Q_{ij}=\frac{\sqrt{r_ir_j}}{r_i+r_j}.
\]
Remark~\ref{rem:cauchy-kernel} gives bounded invertibility of $Q$. The parameter count follows from \eqref{eq:ffc-count-comparison}, and \cref{thm:ffc-critical} proves optimality.
\end{proof}

More generally, the profiles
\begin{equation}\label{eq:ffc-rational-profiles}
 \psi(s)=\frac{s^\alpha}{(1+s)^\beta},\qquad
 \alpha\ge0,\quad\beta>\alpha+\tfrac12,
\end{equation}
with $\psi(0)=1$ for $\alpha=0$ and $\psi(0)=0$ otherwise, are admissible in \cref{thm:ffc-critical}. To see this, we put $a=\alpha+1/2$, $b=\beta-a>0$, and $g_a(u)=e^{au-e^u}$. Positive row and column factors in \eqref{eq:strict-TP} show that $g_a$ and its reflection $\widetilde g_b(u)=g_b(-u)$ are totally positive. Their convolution satisfies
\[
 (g_a*\widetilde g_b)(u)
 =\Gamma(\beta)\frac{e^{au}}{(1+e^u)^\beta}.
\]
Here $\Gamma$ denotes Euler's gamma function. Closure under convolution \cite{Schoenberg51} proves the assertion. These profiles give the vectors $\sqrt t\,(tB)^\alpha(I+tB)^{-\beta}f$, including the resolvent powers $\sqrt t\,(I+tB)^{-\beta}f$ when $\alpha=0$.

The profile $\psi(s)=e^{-s}$ gives $\sqrt t\,e^{-t(I-A)}f$ and kernel $g$. On $(\ker A)^\perp$, replacing $B=I-A$ by $-\log A$ gives the power observations $\sqrt t\,A^t f$ (with $f$ projected onto that subspace). Thus both operator families lead to the same kernel. The Gaussian kernel gives another continuous profile,
$\psi(s)=s^{-1/2}e^{-(\log s)^2}$, with $\psi(0)=0$.

\subsubsection{Integer parameters.}
To obtain integer parameters for the bases above, we round the real parameters and complete the resulting family by finitely many observations.

\begin{proposition}\label{prop:ffc-integer}
Let $(t_k)_{k\ge0}$ be the parameters supplied by
\cref{thm:ffc-critical} or \cref{prop:ffc-resolvent},
with $\psi(s)=(1+s)^{-1}$ in the latter case. Set $n_k=\lceil t_k\rceil$, $T_{\rm int}=\{n_k:k\ge0\}$, and $w_n=n\,\#\{k:n_k=n\}$.
Then the following assertions hold.
\begin{enumerate}[label=\textup{(\roman*)},leftmargin=2em]
\item
The family
$(\sqrt{w_n}\psi(nB)f)_{n\in T_{\rm int}}$ has a Fredholm
analysis operator, and
 $N_{T_{\rm int}}(R)=n_A(R)+O(1)$.
A Riesz basis can be obtained by adjoining finitely many
vectors of the form $A^nf$ and deleting a suitable finite
subset. These modifications change the parameter count
by only $O(1)$.

\item
If $(\psi(nB)f)_{n\in\Nzero}$ is complete, the completing
vectors can be chosen from
 $\{\psi(nB)f:n\in\Nzero\setminus T_{\rm int}\}$.
There is then an integer design within the profile $\psi$
whose weighted sampled family is a Riesz basis and whose
count is $n_A(R)+O(1)$. The completeness assumption holds
for every profile in \eqref{eq:ffc-rational-profiles}.
\end{enumerate}
\end{proposition}
\begin{proof}
We first prove (i). We write $u_k=\log t_k$ and $\eta_k=\log n_k-u_k$. Since $t_k\to\infty$ and $0\le\eta_k\le\log(1+t_k^{-1})$, these shifts are bounded and tend to zero. Before repeated parameters are merged, the difference between the rounded and original analysis matrices has entries
\[
 E_{kj}=\bigl(\varphi(u_k+\eta_k-y_j)-\varphi(u_k-y_j)\bigr)\overline{d_j}.
\]
We choose $D>0$ such that $|u_k-y_k|+\eta_k\le D$ for all $k$, and fix $C_\varphi,\kappa>0$ such that
\[
 |\varphi(u)|\le C_\varphi e^{-\kappa|u|},\qquad u\in\R.
\]
Since $|y_k-y_j|\ge h_1|k-j|$, we obtain
\[
 |E_{kj}|\le2C_\varphi d_+e^{\kappa D}e^{-\kappa h_1|k-j|}.
\]
For each fixed integer $p$, uniform continuity gives
$|E_{\ell+p,\ell}|\le d_+\norm{\varphi(\cdot+\eta_{\ell+p})-\varphi}_\infty\to0$.
Lemma~\ref{lem:band-compact} therefore makes $E$ compact. The rounded analysis is a compact perturbation of an invertible operator and is Fredholm.

If $u_{k+1}-u_k\ge\delta>0$, then
$t_{k+1}-t_k\ge(e^\delta-1)t_k\to\infty$. Rounding creates only finitely many repeated parameters and changes the count by $O(1)$. At a parameter repeated $m_n$ times, a unitary change of the corresponding output coordinates replaces the identical rows by one row with weight $m_nn$ and $m_n-1$ zero rows. Deleting these finitely many zero rows preserves Fredholmness and gives the stated weights.

Original observability and \cref{lem:completion} provide finitely many monomial vectors that make the analysis bounded below. Adding finitely many scalar rows preserves Fredholmness, so the finite-deletion argument following that lemma gives a Riesz basis. This proves (i).

For (ii), we denote by $U$ the analysis operator of the rounded family $(\sqrt{w_n}\psi(nB)f)_{n\in T_{\rm int}}$. Every vector in $\ker U$ is orthogonal to the retained profile vectors. Completeness of the full integer profile family therefore gives
\[
 \ker U\cap\bigcap_{n\in\Nzero\setminus T_{\rm int}}
               \{\psi(nB)f\}^{\perp}=\{0\}.
\]
Lemma~\ref{lem:completion} supplies completing vectors at new parameters, and the same deletion argument stays within this family.

It remains to verify completeness for \eqref{eq:ffc-rational-profiles}. If $x$ is orthogonal to all positive-integer profile vectors, we set $\xi_j=\ip x{e_j}$ and
\[
 H(z)=\sum_j\xi_j\overline{f_j}\,r_j^\alpha(1+zr_j)^{-\beta}.
\]
The function $H$ encodes the orthogonality conditions
as zeros at the positive integers. Its vanishing will
give moment identities for a finite measure supported
on the points $r_j$, allowing us to show that every
coordinate of $x$ is zero.
Principal powers define a bounded holomorphic function on $\Re z>-1/2$: on that half-plane $\Re(1+zr_j)>1/2$, while
$\sum_j|\xi_jf_j|r_j^\alpha<\infty$. Orthogonality gives $H(n)=0$ for every positive integer $n$. These zeros violate the half-plane Blaschke condition \cite[Chapter~2]{Duren70}, so $H\equiv0$.
Termwise differentiation at zero is justified by absolute uniform convergence on compact subsets and yields
\[
 \sum_j\xi_j\overline{f_j}r_j^{\alpha+m}=0,\qquad m\ge0.
\]
The finite complex measure $\sum_j\xi_j\overline{f_j}r_j^\alpha\delta_{r_j}$ annihilates every polynomial on $[0,1]$. Polynomial density makes it zero. Since the $r_j$ are distinct and every $f_jr_j^\alpha$ is nonzero, $x=0$. This proves completeness and (ii).
\end{proof}

Thus, the  profiles in \eqref{eq:ffc-rational-profiles} give integer Riesz bases with critical count entirely within the same function family. For other profiles, the construction may require finitely many monomial observations. Rounding and completion may change the frame bounds.

\subsubsection{Perturbations of the functions.}
The bases in \cref{thm:ffc-critical} also admit square-summable changes of their defining functions, provided completeness is retained. For comparison, Krishtal, Mashreghi, and Miller \cite[Theorem~1.4]{FFC26} treat perturbations of sampled powers on separated spectra in a Stolz domain and an angular sector. Their hypotheses bound the number of sampling exponents in consecutive intervals of a fixed length above and below, bound the $\ell^2$ error uniformly over spectral points, and require the perturbed functions to span a dense subspace of $C(\Gamma)$ for a curve $\Gamma$ containing the spectrum. Here the reference family is the Riesz basis constructed in \cref{thm:ffc-critical}, and the criterion is weighted square-summability together with completeness.

\begin{corollary}\label{cor:ffc-perturbation}
Let $(F_{t_k}(A)f)_{k\ge0}$ be a Riesz basis supplied by
\cref{thm:ffc-critical}, with $F_t$ defined in
\eqref{eq:ffc-profile}.
Let $H_k\in C(\sigma(A))$, and suppose that
\begin{equation}\label{eq:ffc-hs-perturbation}
 \sum_j|f_j|^2
 \sum_k|H_k(\lambda_j)-F_{t_k}(\lambda_j)|^2<\infty.
\end{equation}
If $(H_k(A)f)_k$ is complete, then it is a Riesz basis.
In particular, completeness follows if
$\operatorname{span}\{H_k:k\ge0\}$ is dense in
$C(\sigma(A))$ in the uniform norm.
The simpler bound
\[
 \sup_j\sum_k|H_k(\lambda_j)-F_{t_k}(\lambda_j)|^2<\infty
\]
implies \eqref{eq:ffc-hs-perturbation}.
\end{corollary}
\begin{proof}
Tonelli's theorem and spectral calculus identify the expression in \eqref{eq:ffc-hs-perturbation} with
$\sum_k\|(H_k(A)-F_{t_k}(A))f\|^2$. The two synthesis operators therefore differ by a Hilbert--Schmidt operator. Since the original synthesis operator is invertible, the perturbed one is Fredholm of index zero. Completeness makes its closed range all of $\HH$, and index zero then gives injectivity.

For the density assertion, each $\lambda_j$ is isolated in $\sigma(A)$, so its indicator belongs to $C(\sigma(A))$. Uniform approximation by the span of the $H_k$ shows that the closed span of the vectors $H_k(A)f$ contains $E_A(\{\lambda_j\})f=f_je_j$. Since every $f_j$ is nonzero, the family is complete. The last implication follows from $\sum_j|f_j|^2<\infty$.
\end{proof}

\subsubsection{Common designs for several observations.}
The scalar construction also extends to several observations by tensoring with the identity on the output space. This gives one design for an entire spectral class, with no restriction on the output dimension.

\begin{corollary}\label{cor:ffc-common}
Let $\Lambda=\{\lambda_j:j\ge0\}\subset[0,1)$, where
$\lambda_j$ increases strictly to one. Set
$r_j=1-\lambda_j$, and suppose that
 $\inf_j\log(r_j/r_{j+1})>0$.
Let $\psi\in C([0,\infty);\R)$ be such that
$u\mapsto e^{u/2}\psi(e^u)$ is nonzero, integrable, and
totally positive. Consider the class of observable
positive pairs $(A,C)$ whose spectral measures are
carried by $\Lambda$. The observation spaces and
eigenspaces may be infinite-dimensional.

Then the following assertions hold.
\begin{enumerate}[label=\textup{(\roman*)},leftmargin=2em]
\item
There exists a common real design with weights $w_t=t$
such that
 $\bigl(\sqrt t\,C\psi(t(I-A))\bigr)_{t\in T}$
satisfies frame inequalities for every pair in the class,
and
 $N_T(R)=\#\{j:r_j\ge R^{-1}\}+O(1)$.
\item
There exists an integer design with the same count whose
sampled profile observations satisfy frame inequalities
after finitely many observations
$CA^{\ell_1},\ldots,CA^{\ell_q}$ are adjoined. The parameters,
weights, and exponents can be chosen independently of
$(A,C)$. For every profile in
\eqref{eq:ffc-rational-profiles}, the common integer design
can be chosen entirely from observations
$C\psi(n(I-A))$.

\item
Every common real or integer design for the profile
$\psi$ satisfies
\[
 N_\mu(R)\ge\#\{j:r_j\ge R^{-1}\}-O_\mu(1).
\]
\end{enumerate}
\end{corollary}
\begin{proof}
For each pair, we put $P_j=E_A(\{\lambda_j\})$. The support assumption gives $\sum_jP_j=I$ strongly. Testing original observability on $P_j\HH$ gives
\[
 a(1-\lambda_j^2)\norm{P_jx}^2
 \le\norm{CP_jx}^2
 \le b(1-\lambda_j^2)\norm{P_jx}^2.
\]
We define $Zx=(r_j^{-1/2}CP_jx)_j$. Since
$(1-\lambda_j^2)/r_j=1+\lambda_j\in[1,2)$, orthogonality gives
$a\norm x^2\le\norm{Zx}^2\le2b\norm x^2$.
The map $Z$ retains the normalized observation of each
spectral component separately. The scalar interpolation
matrix then combines these components in the same way
in every output coordinate, which accounts for the
tensor product below.

For (i), we apply \cref{thm:separated-centers} to $\varphi$ and $y_j=\log(1/r_j)$. We denote the selected nodes by $u_k$, put $t_k=e^{u_k}$, and write $M=(\varphi(u_k-y_j))_{k,j\ge0}$. These choices depend only on $\psi$ and $\Lambda$. Spectral expansion yields
\begin{align*}
 ((M\otimes I_\YY)Zx)_k
 &=\sum_j\varphi(u_k-y_j)r_j^{-1/2}CP_jx\\
 &=\sqrt{t_k}\sum_j\psi(t_kr_j)CP_jx
 =\sqrt{t_k}\,C\psi(t_k(I-A))x.
\end{align*}
The series converges in $\YY$ by continuous functional calculus and boundedness of $C$. Thus the sampled analysis is $(M\otimes I_\YY)Z$. The first factor is boundedly invertible, and the second is bounded and bounded below. This proves the frame inequalities, also for infinite-dimensional $\YY$. Bounded displacement and separation give the stated count.

For (ii), we use the canonical scalar pair on $\ell^2(\Nzero)$,
\[
 A_0e_j=\lambda_je_j,\qquad
 f^{(0)}=\sum_j\sqrt{r_j}e_j,\qquad C_0x=\ip x{f^{(0)}}.
\]
With $h=\inf_{j\ge0}\log(r_j/r_{j+1})>0$, we have
$r_j\le r_0e^{-hj}$, so $f^{(0)}$ is well defined.
Moreover, $r_{j+1}/r_j\le e^{-h}<1$, and hence
$(\lambda_j)$ satisfies the Carleson condition by
\cite[Proposition~1.2(i)]{CHPS24}. Since
\[
 \sqrt{r_j}
 =(1+\lambda_j)^{-1/2}\sqrt{1-\lambda_j^2},
 \qquad
 2^{-1/2}\le(1+\lambda_j)^{-1/2}\le1,
\]
the canonical pair is observable by
\cite[Theorem~1.1]{CHPS24}.

We apply \cref{prop:ffc-integer} to this canonical scalar pair and denote the resulting Riesz-basis analysis operator
by $H$. The finite completion and deletion are performed before tensoring, since a nonzero kernel or cokernel of the rounded scalar analysis operator would become
infinite dimensional upon tensoring with an
infinite-dimensional output space. The profile rows of $H$ and any retained monomial rows
have the forms $(\sqrt{w_nr_j}\,\psi(nr_j))_j$ and
$(\sqrt{r_j}\lambda_j^\ell)_j$, respectively.
After tensoring with $I_\YY$ and composing with $Z$,
these rows become, respectively,
$\sqrt{w_n}C\psi(n(I-A))$ and $CA^\ell$. The resulting analysis operator is
$(H\otimes I_\YY)Z$, which is bounded and bounded below. All choices were made on the canonical pair, so they are independent of $(A,C)$. The count follows from \cref{prop:ffc-integer} and \eqref{eq:ffc-count-comparison} for $A_0$.
For the profiles in \eqref{eq:ffc-rational-profiles}, the same proposition permits all completing vectors to be chosen within the profile family, proving the final part of (ii).

Finally, every common design is admissible for the canonical pair. The lower bound in \cref{thm:ffc-critical} and \eqref{eq:ffc-count-comparison} give
\[
 N_\mu(R)\ge n_{A_0}(R)-O_\mu(1)
 =\#\{j:r_j\ge R^{-1}\}-O_\mu(1).
\]
This proves (iii).
\end{proof}

These designs are optimal for the entire spectral class. For a fixed pair, the count can be smaller. Our forthcoming research will study how the relative positions of the observation ranges determine the optimum and its attainment.

\section{A normal orbit with infinite counting coefficient}\label{sec:normal-obstruction}
The scalar result in \cref{thm:critical-single} uses positivity to obtain
separation of the logarithmic decay rates. For normal operators, angular
separation can yield an observable orbit even when consecutive decay
rates have ratio tending to one. The following example shows that the
radial spectral count can then underestimate the number of observations
required for stable recovery by an unbounded factor.

In this section, $A$ is a bounded normal operator and
$|A|=(A^*A)^{1/2}$. We use the radial count
\[
 n_A(R)=\rank E_{|A|}((0,e^{-1/R}]),\qquad R>0,
\]
which agrees with \eqref{eq:spectral-count} when $A$ is positive.
Real powers are defined by Borel spectral calculus using
$z^t=|z|^t e^{it\arg z}$ for $z\ne0$, with
$\arg z\in(-\pi,\pi]$, $0^t=0$ for $t>0$, and $A^0=I$.
For $Cx=\ip{x}{f}$, a design $\mu=\sum_{t\in T}w_t\delta_t$
is admissible if
\begin{equation}\label{eq:normal-admissibility}
 \alpha\norm x^2
 \le\sum_{t\in T}w_t\left|\ip{A^tx}{f}\right|^2
 \le\beta\norm x^2,\qquad x\in\HH,
\end{equation}
for some $0<\alpha\le\beta<\infty$. Thus, the sampled frame consists
of the vectors $\sqrt{w_t}(A^t)^*f$.

\begin{proposition}\label{prop:normal-infinite}
There exist a normal operator $A$ with $\norm A=1$ on an
infinite-dimensional Hilbert space $\HH$ and a vector $f\in\HH$
such that both $(A^nf)_{n\ge0}$ and $((A^*)^nf)_{n\ge0}$ are frames,
but every admissible real design for $Cx=\ip{x}{f}$ satisfies
\[
 \frac{N_\mu(R)}{n_A(R)}\longrightarrow+\infty
 \qquad(R\to\infty).
\]
More precisely, $n_A(R)\asymp R/(\log R)^2$, whereas
$N_\mu(R)\ge c_\mu R/\log R$ for all sufficiently large $R$,
with $c_\mu>0$.
\end{proposition}
\begin{proof}
We fix $\varepsilon=10^{-3}$ and, for $j\ge2$, put
\[
 a_j=\frac{\varepsilon}{j(\log j)^2},\qquad
 \theta_j=4\sum_{k=j}^\infty a_k,\qquad
 \lambda_j=(1-a_j)e^{i\theta_j},\qquad
 q_j=1-|\lambda_j|^2.
\]
Our choice of
$\varepsilon$ gives $a_j<1/4$ and $0<\theta_j<\pi/2$.
On an orthonormal basis $(e_j)_{j\ge2}$, we define
$Ae_j=\lambda_je_j$ and $f=\sum\limits_{j\ge2}\sqrt{q_j}e_j$
(the summability of $(a_j)$ gives $f\in\HH$). $A$ is normal
and $\|A\|=1$.

We first verify the integer frame property. The boundary arcs
 $J_j=\{e^{i\theta}:|\theta-\theta_j|<a_j/2\}$
are pairwise disjoint, since
 $\theta_j-\theta_{j+1}=4a_j>
 \frac{a_j+a_{j+1}}{2}$.
Their lengths are $a_j=1-|\lambda_j|$. Hence the
disjoint-arc criterion of Stray and \O{}yma
\cite[pp.~77--78, proof of Theorem~1]{SO78}
shows that $(\lambda_j)$ is an interpolating sequence
for $H^\infty(\mathbb D)$,
\(\mathbb D=\{z\in\mathbb C: |z| < 1\}\), or equivalently a Carleson sequence.
Since $q_j=1-|\lambda_j|^2$,
\cite[Theorem~2.2]{CMPP20} implies that $(A^nf)_{n\ge0}$
is a frame. Conjugation preserves the Carleson condition,
so $((A^*)^nf)_{n\ge0}$ is also a frame.
In particular, $(A,C)$ is observable.

We set $s_j=-\log(1-a_j)$. Then $s_j\sim a_j$ and
$s_j/s_{j+1}\to1$. Inverting $j(\log j)^2\sim\varepsilon R$ gives
\begin{equation}\label{eq:normal-spectral-count}
 n_A(R)=\#\{j:s_j\ge R^{-1}\}
 \sim\frac{\varepsilon R}{(\log R)^2}.
\end{equation}
The function
\[
 P(t)=\norm{(A^t)^*f}^2
     =\sum_{j\ge2}q_je^{-2s_jt}
\]
decreases with $t$. For $j\ge\lceil R\rceil$ and large $R$,
the quantities $Rs_j$ are uniformly bounded. Since $q_j\ge a_j$,
the integral estimate for
$\sum_{j\ge\lceil R\rceil}1/(j(\log j)^2)$ gives
\begin{equation}\label{eq:normal-orbit-tail}
 P(R)\ge\frac{c\varepsilon}{\log R}
\end{equation}
for some $c>0$ and all sufficiently large $R$.

We now fix an admissible design with bounds $\alpha,\beta$
in \eqref{eq:normal-admissibility}, and put
\[
 W_\mu(R)=\sum_{\substack{t\in T\\t\le R}}w_t.
\]
Testing the frame inequalities on $e_j$ gives
 $\alpha\le q_j\sum_{t\in T}w_te^{-2s_jt}\le\beta$.
Since $q_j=1-e^{-2s_j}$, we have
$s_j\le q_j\le2s_j$ for all sufficiently large $j$.

We first obtain an upper bound for $W_\mu$ by applying the Bessel inequality to a single eigenvector, as in the
proof of \cref{cor:riesz-count}.
Since $s_j/s_{j+1}\to1$, for every sufficiently large $R$
there is an index $k$ such that
$1/(2R)\le s_k\le 1/R.$
Consequently,
\[
 \beta
 \ge q_k\sum_{\substack{t\in T\\t\le R}}w_te^{-2s_kt}
 \ge \frac{e^{-2}}{2R}W_\mu(R),
\]
and hence $W_\mu(R)\le2e^2\beta R$.

For the lower bound, we fix $b>0$ so that
$8\beta e^{-b}\le\alpha/2$.
The same ratio condition allows us, for every sufficiently
large $R$, to choose indices $j$ and $\ell$ satisfying
$b/R\le s_j\le 2b/R,\ s_j/4\le s_\ell\le s_j/2.$
Thus $q_j/q_\ell\le8$ and
$2(s_j-s_\ell)R\ge b$. The Bessel estimate on $e_\ell$
bounds the late observations of $e_j$:
\[
 \begin{aligned}
 q_j\sum_{\substack{t\in T\\t>R}}w_te^{-2s_jt}
 \le
 \frac{q_j}{q_\ell}e^{-2(s_j-s_\ell)R}
 q_\ell\sum_{\substack{t\in T\\t>R}}w_te^{-2s_\ell t}
 \le8\beta e^{-b}\le\frac{\alpha}{2}.
 \end{aligned}
\]
The lower frame bound on $e_j$ therefore yields
$\alpha\le q_jW_\mu(R)+\alpha/2,\ W_\mu(R)\ge\alpha/(2q_j)\ge\alpha R/(8b).$
This proves $W_\mu(R)\asymp R$.

Each individual observation in \eqref{eq:normal-admissibility}
has squared norm at most $\beta$, so $w_tP(t)\le\beta$.
Since $P$ decreases,
\[
 W_\mu(R)\le\frac{\beta N_\mu(R)}{P(R)}.
\]
Combining this inequality with \eqref{eq:normal-orbit-tail}
and the lower bound for $W_\mu$ gives
\[
 N_\mu(R)\ge\frac{P(R)W_\mu(R)}{\beta}
 \ge c_\mu\frac{R}{\log R}.
\]
Together with \eqref{eq:normal-spectral-count}, this proves
the assertion.
\end{proof}

The unweighted integer design is admissible in this example, so
both the real and integer admissible classes are nonempty.
Nevertheless, every admissible design has an infinite counting
coefficient relative to the radial count. Thus, positivity is
essential even for the existence of a design with
$N_\mu(R)=O(n_A(R))$ in the scalar problem.

\section*{Declaration of AI use}
The authors used ChatGPT (OpenAI) to assist with developing and checking mathematical arguments, identifying relevant
literature, and drafting the manuscript.
The authors verified the mathematical arguments and take full responsibility for the results and the final
content of the paper.

\raggedbottom
\renewcommand{\bibliofont}{\small}

\end{document}